\documentclass{article}
\usepackage{amsmath}
\usepackage{amsfonts,amssymb,mathrsfs}
\usepackage{theorem}
\usepackage{tikz}
\usepgflibrary{arrows}

\renewcommand{\Im}{\mathrm{Im}}

\input amssym.def

\numberwithin{equation}{section}

\numberwithin{figure}{section}

\makeatletter
\newcommand\qedsymbol{\hbox{$\Box$}}
\newcommand\qed{\relax\ifmmode\Box\else
  {\unskip\nobreak\hfil\penalty50\hskip1em\null\nobreak\hfil\qedsymbol
  \parfillskip=\z@\finalhyphendemerits=0\endgraf}\fi}

\newenvironment{proof}[1][{}]{\par\noindent Proof{#1}. }{\qed}

\newcommand{\bfzero}{{\mathbf{0} }}

\newcommand{\PV}{\mathsf{PV}}

\newcommand{\dGra}{{\mathsf{dGra} }}
\newcommand{\KGra}{{\mathsf {KGra} }}

\newcommand{\Cobar}{{\mathrm{Cobar} }}

\renewcommand{\c}{{\circ}}
\renewcommand{\i}{\sqrt{-1}\,}

\newcommand{\dfGC}{{\mathsf{dfGC}}}

\newcommand{\dgra}{\mathrm{dgra}}

\newcommand{\Lie}{{\mathsf{Lie}}}

\newcommand{\Conf}{{\mathsf{Conf}}}

\newcommand{\OC}{{\mathsf{OC}}}
\newcommand{\br}{{\mathrm{br}}}

\newcommand{\Ch}{{\mathsf{Ch}}}
\newcommand{\grVect}{{\mathsf{grVect}}}

\newcommand{\A}{\mathscr{A}}
\newcommand{\B}{\mathscr{B}}
\newcommand{\C}{\mathscr{C}}

\newcommand{\Op}{{\mathbb{OP}}}

\newcommand{\oc}{{\mathfrak{o} \mathfrak{c}}}

\newcommand{\Conv}{{\mathrm{Conv}}}

\newcommand{\Sh}{\mathrm{Sh}}

\newcommand{\End}{{\mathsf{E n d}}}
\newcommand{\Hom}{{\mathrm{Hom}}}
\newcommand{\Arg}{{\mathrm{Arg}}}

\newcommand{\Alt}{{\mathrm{Alt}}}

\newcommand{\inv}{{\mathrm {inv} }}
\newcommand{\Hoch}{\mathrm{Hoch}}

\newcommand{\id}{\mathrm{id}}
\newcommand{\ad}{\mathrm{ad}}

\newcommand{\tm}{\tilde{m}}

\newcommand{\ti}[1]{{\tilde{#1}}}

\newcommand{\sh}{\sharp}

\newcommand{\dia}{\diamond}

\newcommand{\tbeta}{\tilde{\beta}}

\newcommand{\Cbu}{C^{\bullet}}

\newcommand{\al}{{\alpha}}
\newcommand{\la}{{\lambda}}
\newcommand{\vr}{{\varrho}}

\newcommand{\bul}{{\bullet}}
\newcommand{\hs}{{\heartsuit}}

\newcommand{\ww}{{\circ\, \circ}}
\newcommand{\ed}{{\bullet\hspace{-0.05cm}-\hspace{-0.05cm}\bullet}}

\newcommand{\md}{{\mathfrak{d}}}

\newcommand{\mc}{{\mathfrak{c}}}
\newcommand{\mo}{{\mathfrak{o}}}

\newcommand{\bC}{\overline{C}}

\newcommand{\si}{{\sigma}}
\newcommand{\ga}{{\gamma}}
\newcommand{\vf}{{\varphi}}
\newcommand{\ve}{{\varepsilon}}

\newcommand{\ka}{{\kappa}}
\newcommand{\G}{{\Gamma}}
\newcommand{\cF}{{\cal F}}
\newcommand{\pa}{{\partial}}

\newcommand{\st}{{\mathsf{t}}}

\newcommand{\bs}{{\mathbf{s}}}
\newcommand{\bsi}{{\mathbf{s}^{-1}\,}}

\newcommand{\bq}{{\mathbf{q}}}

\newcommand{\cC}{{\cal C}}

\newcommand{\cL}{{\mathcal{L}}}
\newcommand{\cD}{\mathcal{D}}
\newcommand{\cA}{{\cal A}}

\newcommand{\cR}{\mathcal{R}}

\newcommand{\cO}{{\cal O}}
\newcommand{\cV}{{\cal V}}

\newcommand{\bz}{{\bar z}}

\newcommand{\bbK}{{\mathbb K}}

\newcommand{\bbC}{{\mathbb C}}
\newcommand{\bbR}{{\mathbb R}}
\newcommand{\bbZ}{{\mathbb Z}}
\newcommand{\bbQ}{{\mathbb Q}}

\newcommand{\La}{{\Lambda}}

\newcommand{\n}{{\nabla}}
\newcommand{\te}{\theta}

\newcommand{\de}{{\delta}}

\newcommand{\sgn}{{\rm s g n}}

\date{}

\newtheorem{lem}{Lemma}[section]

\newtheorem{defi}[lem]{Definition}

\newtheorem{thm}[lem]{Theorem}
\newtheorem{cor}[lem]{Corollary}

\newtheorem{prop}[lem]{Proposition}
\newtheorem{claim}[lem]{Claim}

\theorembodyfont{\rm}
\newtheorem{example}[lem]{Example}
\newtheorem{pty}[lem]{Property}
\newtheorem{remark}[lem]{Remark}

\title{A Formality Quasi-isomorphism for Hochschild Cochains 
over Rationals Can Be Constructed Recursively}

\author{Vasily Dolgushev}
\date{}

\begin{document}

\large

\maketitle

\begin{center}
{\it Belatedly to Volodya Rubtsov on the occasion of his 60th birthday.}
\end{center}

~\\[-0.3cm]

\begin{abstract}
It is believed  \cite{FW-irrat}, \cite[Section 4.1]{K-mot} that, among the coefficients entering 
Kontsevich's formality quasi-isomorphism \cite{K}, there are irrational (possibly even transcendental)
numbers. In this paper, we prove that a formality quasi-isomorphism for Hochschild cochains 
of a polynomial algebra over $\bbQ$ can be constructed recursively. The proof that the 
proposed recursive algorithm works, is based on the existence of formality quasi-isomorphism 
over $\bbR$.  However, the algorithm requires no explicit knowledge of the coefficients 
entering Kontsevich's construction. Although this algorithm completely bypasses 
Tamarkin's approach \cite{Hinich}, \cite{Dima-Proof}, the construction is inspired 
by Proposition 5.8 from classical paper \cite{Drinfeld} by V. Drinfeld. 
\end{abstract}

\tableofcontents

\section{Introduction}
It is notoriously hard \cite{AM}, \cite{Sahi}, \cite{FW-irrat}, \cite{Kathotia},  \cite[Section 4.1]{K-mot}, 
\cite{Polyak},  \cite{Sh}, \cite{VdB} to compute the coefficients (weights) entering Kontsevich's formality 
quasi-isomorphism for Hochschild cochains \cite{K}. This difficulty is a stumbling block  
for direct applications of Kontsevich's construction to concrete questions in deformation quantization.   
In addition, it is very likely \cite{FW-irrat}, \cite[Section 4.1]{K-mot} that some of these weights are 
irrational and possibly even transcendental numbers. 

The goal of this paper is to propose a construction which, in some sense, demystifies formality 
quasi-isomorphisms for Hochschild cochains and hence unlocks new tools for applications 
of formality quasi-isomorphisms for Hochschild cochains to deformation quantization.   
More precisely, we propose an algorithm whose output is an infinite sequence of 
approximations which ``converge'' to a formality quasi-isomorphism for Hochschild cochains 
(defined over $\bbQ$). Given any such approximation, the construction of 
the next (better) approximation boils down to solving a finite dimensional linear system.
The proof that this algorithm works, {\it does use} the existence of formality quasi-isomorphism 
over $\bbR$ (established in \cite{K}). However, the algorithm requires {\it no explicit knowledge} of 
these mysterious weights entering Kontsevich's construction.  

Recall that Tamarkin's construction \cite[Section 2]{DP}, \cite{Hinich}, \cite{Dima-Proof} gives us a map 
from the set of Drinfeld associators to the set of homotopy classes of formality quasi-isomorphisms 
for Hochschild cochains of the polynomial algebra. In this respect, the proposed construction is inspired by 
Proposition 5.8 from classical paper\footnote{See also Theorem 4 and Corollary 4.1 in D. Bar-Natan's beautiful paper \cite{Bar-Natan}.}
\cite{Drinfeld} by V. Drinfeld.  The proof that the proposed algorithm works, is also based on the action of the 
full directed graph complex $\dfGC$ \cite{stable}, \cite{dfGC}, \cite[Section 5]{K-conj}, 
\cite{Thomas} on the formality quasi-isomorphisms 
for Hochschild cochains.

\subsection{A colloquial description of the main result}

The main ingredient of the construction, an $m$-th approximation to a formality quasi-isomorphism,  
can be described colloquially without using the language of operads. 

To give this informal description, we denote by $\PV$ the vector space 
of polyvector fields on an affine space and recall that the problem of 
constructing a formality quasi-isomorphism for Hochschild cochains 
$\Cbu(A) : = \Cbu(A,A)$ 
of a polynomial algebra $A$ is equivalent to constructing 
a collection of maps
\begin{equation}
\label{U-n-k-intro}
U_{n,k} ~:~ \PV^{\otimes\, n} \otimes A^{\otimes k} ~\to~ A
\end{equation}
such that 
\begin{itemize}

\item the collection $\{U_{1,k}\}_{k \ge 0}$ assembles into the standard Hochschild-Kostant-Rosenberg 
embedding $\PV \hookrightarrow \Cbu(A)$, 

\item for every $n \ge 2, ~ k \ge 0$, the map $U_{n,k}$ is compatible with the action of $S_n$ on 
$\PV^{\otimes\, n}$ and 

\item the collection of maps \eqref{U-n-k-intro} satisfies a sequence of certain quadratic relations
\begin{equation}
\label{cR-n-k}
\cR_{n,k} \big( \{U_{n_1,k_1}\}_{n_1 \ge 1, ~k_1 \ge 0} \big) = 0\,, \qquad n \ge 1, ~~ k \ge 0. 
\end{equation}

\end{itemize}

The right hand side of \eqref{cR-n-k} is an element of 
$$
\Hom(\PV^{\otimes\, n} \otimes A^{\otimes\, k}\,,\, A) 
$$
which is expressed in terms of maps  $\{U_{n_1,k_1}\}_{n_1 \ge 1, ~k_1 \ge 0}$ 
using the obvious compositions of Hom's, the multiplication on $A$, and the Schouten 
bracket on $\PV$.

For an integer $m \ge 2$, we introduce the following 
subsets of\footnote{Examples of $\A_m$ and $\B_m$, for $m=5$, are depicted in 
figure \ref{fig:areas-m-approx}.} 
$\bbZ_{\ge 1} \times \bbZ_{\ge 0}$: 
\begin{equation}
\label{A-m-intro}
\A_m : = \big\{ (n,k) \in \bbZ_{\ge 1} \times \bbZ_{\ge 0} ~ \big|~ 2n + k \le 2 m + 1   \big\} \cup
\big\{ (1,k)~ \big|~ \textrm{for any } k \ge 0  \big\},
\end{equation}
\begin{equation}
\label{B-m-intro}
\B_m : = \{(m+1, 0),\, (m, 2),\,  (m-1, 4),\, \dots,\, (2, 2m-2) \}.
\end{equation}

We observe that, for every $(n,k) \in \A_m \cup \B_m$, the relation $\cR_{n,k}$ 
depends only on the maps $\{ U_{n_1, k_1} \}_{(n_1, k_1) \in \A_m}$. 
Due to this observation, we can define (colloquially) an $m$-th approximation $U^{(m)}$
to a formality quasi-isomorphism as a family of maps (defined over $\bbQ$)
$$
U^{(m)}_{n_1,k_1} ~:~ \PV^{\otimes\, n_1} \otimes A^{\otimes k_1} ~\to~ A
$$
such that 
\begin{itemize}

\item $U^{(m)}_{n_1,k_1} = 0$ for all $(n_1, k_1) \notin \A_m$,

\item the relations
$$
\cR_{n,k} \big( \{U^{(m)}_{n_1,k_1} \}_{n_1\ge 1,~ k_1 \ge 0} \big) = 0
$$ 
are satisfied for all $(n,k) \in \A_m \cup \B_m$, and 

\item the maps $\{ U^{(m)}_{1,k_1}\, ,\,  U^{(m)}_{2,k_1} \}_{k_1 \ge 0}$
satisfy some technical conditions. 

\end{itemize}

The main result of this paper (see Theorem \ref{thm:main}) can be stated (colloquially) as follows: 

~\\
{\it The second approximation exists. Given an $m$-th approximation $U^{(m)}$ to a formality 
quasi-isomorphism, one can construct an $(m+1)$-th approximation $U^{(m+1)}$ such that
$$
U^{(m+1)}_{n,k} = U^{(m)}_{n,k}\,, \qquad \forall ~ (n,k) \in \A_m\,. 
$$
Finally, to construct $U^{(m+1)}$, one needs to solve a finite dimensional linear system.}\\

Although an $m$-th approximation $U^{(m)}$ is not an $L_{\infty}$ quasi-isomorphism 
from $\PV$ to $\Cbu(A)$, it can still be used to construct associative star products in 
$A[\ve]/(\ve^m)$, where $\ve$ is the formal deformation parameter. 
More precisely\footnote{See Theorem \ref{thm:star-approx} for a more general statement.},
\begin{thm}
\label{thm:star-prod}
For every Poisson structure $\ka \in \PV^2$, the formula 
\begin{equation}
\label{star-intro}
a * b : = a b + \sum_{n = 1}^{m-1} \frac{\ve^n}{n!} U^{(m)}_{n,2}(\underbrace{\ka, \ka, \dots, \ka}_{n~\emph{times}} ; a,b), 
\qquad a,b \in A
\end{equation}
defines an associative multiplication on $A[\ve]/(\ve^m)$. Moreover, \eqref{star-intro} is 
a truncation modulo $(\ve^m)$ of an honest star product on $A[[\ve]]$. 
\end{thm}

\subsection{Sketch of the proof}
Let us denote by $U^K$ Kontsevich's formality quasi-isomorphism from \cite{K}
defined over $\bbR$. Due to Claim \ref{cl:Konts-2} from Appendix \ref{app:Konts-ptys}, 
$$
U^{K}_{2,0} = 0 \qquad \textrm{and} \qquad  U^{K}_{2,1} = 0
$$
and hence $U^{K}_{2,0}$ and $U^{K}_{2,1}$ are defined over $\bbQ$. 
Thus a second approximation $U^{(2)}$ to a formality quasi-isomorphism exists.  

Let $U^{(m)}$ be an $m$-th approximation to a formality quasi-isomorphism. 
In general, $U^{(m)}_{n,k}  \neq U^{K}_{n,k}$ for $(n,k) \in \A_m$.  
However, there exists a sequence of genuine formality quasi-isomorphisms 
(defined over $\bbR$)
$$
\{ U^{\hs ,  \tm} \}_{1 \le \tm \le m} 
$$ 
such that 
\begin{itemize}

\item $U^{\hs, 1} = U^{K}$ and

\item $U^{\hs ,  \tm}_{n,k}  =  U^{(m)}_{n,k} $ for all $(n,k) \in \A_{\tm}$ and $2 \le \tm \le m$. 

\end{itemize}

$U^{\hs, \tm}$ is obtained from $U^{\hs, \tm-1}$ by using the homotopy equivalence 
and the action of the full directed graph complex \cite{stable}, \cite{dfGC}. 

It is the existence of $U^{\hs, m}$, which guarantees that the linear system 
for constructing the better approximation $U^{(m+1)}$ is consistent over $\bbR$. 
Thus, since both the coefficient matrix and the right hand side of this linear system 
are defined over $\bbQ$, it is also consistent over $\bbQ$. Therefore, $U^{(m)}$ can 
be extended to the next approximation $U^{(m+1)}$.

\subsection{The organization of the paper}
In the last part of the introduction, we go over the notational 
conventions. In Section \ref{sec:OC-KGra-SFQ}, we give a brief 
reminder of the Kajiura-Stasheff operad $\OC$, the operad $\KGra$
and the notion of a stable formality quasi-isomorphism (SFQ) 
following \cite{stable}.  In Section \ref{sec:OC-KGra-SFQ}, we
also recall Kontsevich's construction \cite{K} of the first 
example of an SFQ (over $\bbR$). 

In Section \ref{sec:prep}, we prove Propositions \ref{prop:al-al} and 
\ref{prop:exp-xi-beta}. Proposition \ref{prop:al-al} allows us to introduce 
the notion of an $m$-th approximation to an SFQ and Proposition 
\ref{prop:exp-xi-beta} plays an important role in the proof of the 
main theorem. 

In Section \ref{sec:main-thm}, we introduce the notion of an $m$-th 
approximation to an SFQ (for $m \ge 2$) and prove the existence and 
uniqueness of the second approximation to an SFQ. In this section, we also 
formulate the main result of this paper (see Theorem \ref{thm:main}) and 
deduce it from technical Proposition \ref{prop:main}.

In Section \ref{sec:star}, we prove that an $m$-th approximation to an SFQ 
can be used to construct star products on an affine space modulo $(\ve^m)$, 
where $\ve$ denotes a formal deformation parameter.

Section \ref{sec:proof-main} is devoted to the proof of 
Proposition \ref{prop:main} which, in turn, implies the main theorem 
(Theorem \ref{thm:main}). Finally, Appendix \ref{app:Konts-ptys} is devoted 
to the proofs of some technical properties of Kontsevich's SFQ from \cite{K}.

\subsection{Notation and conventions}
Many notational conventions are borrowed from \cite{stable}.
The base field $\bbK$ has characteristic zero and, in this paper,
$\bbK$ is often $\bbQ$ or $\bbR$.  

The underlying symmetric monoidal category  for our algebraic 
structures is either the category $\grVect_{\bbK}$ of $\bbZ$-graded
$\bbK$-vector spaces or the category 
$\Ch_{\bbK}$ of unbounded cochain complexes of $\bbK$-vector spaces.  

For a cochain complex $\cV$ we denote 
by $\bs \cV$ (resp. by $\bs^{-1} \cV$) the suspension (resp. the 
desuspension) of $\cV$\,. In other words, 
$$
\big(\bs \cV\big)^{\bul} = \cV^{\bul-1}\,,  \qquad
\big(\bs^{-1} \cV\big)^{\bul} = \cV^{\bul+1}\,. 
$$
For a homogeneous vector $v$ in a graded vector space (or a cochain complex) 
$\cV$, the notation $|v|$ is reserved for the degree of $v$.  

$\Cbu(A)$  denotes the Hochschild cochain
complex of an associative algebra  (or more generally an $A_{\infty}$-algebra) $A$
with coefficients in $A$\,. For a commutative ring $R$ and an $R$-module $V$ we denote 
by $S_R(V)$ the symmetric algebra of $V$ over $R$\,. For a vector $\xi$ in a (dg) Lie algebra 
$\cL$, the notation $\ad_{\xi}$ is reserved for the adjoint action of $\xi$, i.e. 
$\ad_{\xi} (\eta) : =  [\xi, \eta]$\,.

Given an operad $\cO$, we denote by $\circ_i$ the elementary 
operadic insertions: 
$$
\circ_i : \cO(n) \otimes \cO(k) \to \cO(n+k-1)\,, \qquad 1 \le i \le n.
$$
For a collection $Q$, the notation $\Op(Q)$ is reserved for the 
free operad generated by this collection.  

The symmetric group on $n$ letters in denoted by $S_n$ and
the notation $\Sh_{p,q}$ is reserved for the set of $(p,q)$-shuffles 
in $S_{p+q}$\,. A graph is {\it directed} if each edge 
carries a chosen direction. A graph $\G$ with 
$n$ vertices is called {\it labeled} if $\G$ is equipped with a bijection between 
the set of its vertices and the set $\{1,2, \dots, n\}$\,. In this paper we 
consider exclusively graphs without loops (i.e. cycles of length one).

We will freely use the conventions for colored 
(co)operads and colored pseudo- (co)operads from \cite[Section 2]{stable}. 
For example, for a (colored) pseudo-cooperad $\cC$ and a (colored) pseudo-operad $\cO$, the 
notation
$$
\Conv(\cC, \cO)
$$
is reserved for the convolution dg Lie algebra \cite[Section 2.3]{stable}, 
\cite[Section 4]{notes}, \cite{MVnado}. 

For most of colored (co)operads and pseudo- (co)operads used in this paper 
the ordinal of colors will have only two element $\mc$ and $\mo$ for 
which we set $\mc < \mo$. Just as in \cite{stable}, solid edges of colored planar trees
are the edges which carry the color $\mc$ and dashed edges are the edges which
carry the color $\mo$. 

For a two colored operad $\cO$, the notation $\cO(n,k)^{\mc}$ (resp. $\cO(n,k)^{\mo}$) is 
reserved for the space of ``operations'' with $n$ inputs of the color $\mc$, $k$ inputs of 
the color $\mo$, and with the output carrying the color $\mc$ (resp. $\mo$).

Using ``arity'' we can equip the convolution Lie algebra $\Conv(\cC, \cO)$ 
with the natural descending filtration 
$$
\Conv(\cC, \cO) =  \cF_{-1}\, \Conv(\cC, \cO) \supset \cF_0\, \Conv(\cC, \cO) \supset \cF_1\, \Conv(\cC, \cO)
\supset  \dots\,,
$$
where
\begin{equation}
\label{Conv-filtr}
\cF_m\, \Conv(\cC, \cO)=
\end{equation} 
$$
 \big\{ f \in \Conv(\cC, \cO) \quad \big| \quad
f \Big|_{\cC(\bq)} \,=\, 0 \quad \forall~~
\textrm{corollas} ~~ \bq \quad \textrm{satisfying} ~ |\bq| \le m \big\}\,,
$$
and $|\bq|$ is the total number of incoming edges of the corolla $\bq$.

This filtration is compatible with the Lie bracket and  $\Conv(\cC, \cO)$ is complete 
with respect to this filtration.  Namely,
\begin{equation}
\label{complete}
\Conv(\cC, \cO) = \lim_{m} \Conv(\cC, \cO) ~\Big/~   \cF_m\, \Conv(\cC, \cO).
\end{equation}

We also introduce an additional descending 
filtration  $\cF^{\chi}_{\bul}$ on the convolution Lie algebra $\Conv(\cC, \cO)$ 
for each color $\chi$:
$$
\Conv(\cC, \cO) =  \cF^{\chi}_{-1}\, \Conv(\cC, \cO) \supset 
 \cF^{\chi}_0\, \Conv(\cC, \cO) \supset \cF^{\chi}_1\, \Conv(\cC, \cO)
 \supset \dots\,,
$$
where
\begin{equation}
\label{Conv-filtr-chi}
\cF^{\chi}_m\, \Conv(\cC, \cO)
\end{equation} 
consists of vectors  $f \in \Conv(\cC, \cO)$ satisfying this condition: 
\begin{equation}
\label{color-in-minus-out}
f \Big|_{\cC(\bq)} \,=\, 0  \quad \textrm{if} \quad \sh^{in}_{\chi}(\bq) - \sh^{out}_{\chi}(\bq) \le m-1. 
\end{equation}
Here $ \sh^{in}_{\chi}(\bq) $ is the number of incoming edges of the corolla $\bq$ 
which carry the color $\chi$ and  $ \sh^{out}_{\chi}(\bq) $ is the number of 
outgoing edges of the corolla $\bq$ which carry the color $\chi$. 
Note that  $ \sh^{out}_{\chi}(\bq) $ is either $1$ or $0$ because every corolla has 
exactly one outgoing edge.

The filtration \eqref{Conv-filtr-chi} is compatible with the Lie bracket on $\Conv(\cC, \cO)$
and  $\Conv(\cC, \cO)$ is complete with respect to 
this filtration.

We denote by $\La$ the endomorphism operad of the $1$-dimensional
vector space $\bs^{-1} \bbK$ placed in degree $-1$
\begin{equation}
\label{La}
\La = \End_{\bs^{-1} \bbK}\,.
\end{equation}
In other words,
$$
\La(n) = \bs^{1-n} \sgn_n\,,
$$
where $\sgn_n$ is the sign representation for the symmetric group $S_n$
We observe that the collection $\La$ is also naturally a cooperad. 

For a dg operad (resp. a dg cooperad) $P$ in we denote by $\La P$ the 
dg operad (resp. the dg cooperad) which is obtained from $P$ via tensoring  
with $\La$, i.e.
\begin{equation}
\label{La-P}
\La P (n) = \bs^{1-n} P(n) \otimes \sgn_n\,.
\end{equation}

\begin{remark}
\label{rem:brack-pa-cobar}
Let $P$ be a pseudo-operad and $\cO$ be an operad. 
Let us recall \cite[Section 2.5]{stable} that the Lie bracket on 
$\Conv(P, \cO)$ can be expressed in terms of the cobar
differential $\cD_{\Cobar}$ on $\Cobar(\cC)$, where $\cC$ is the cooperad 
which is obtained from $P$ via adjoining the counit. 
More precisely, for $f, g \in \Conv(P, \cO)$ and $X \in P$ we have  
\begin{equation}
\label{brack-pa-cobar}
[f,g] (X) =(-1)^{|g|} \mu \big( f \bs^{-1} \otimes  g \bs^{-1} (\cD_{\Cobar}(\bs X)) \big) 
- (-1)^{|f| |g|} (f \leftrightarrow g)\,,
\end{equation}
where $f \bs^{-1}$ and $g \bs^{-1}$ act in the obvious way 
on the tensor factors of $\cD_{\Cobar} (\bs X) \in \Op (\bs\, P)$ 
and $\mu$ denotes the  multiplication map
$$
\mu : \Op(\cO) \to \cO.
$$
\end{remark}

~\\

\noindent
\textbf{Anniversary note:} If you have never met Volodya Rubtsov then I 
strongly suggest you visit him in the City of Angers in France.
You will meet a person who radiates an amazing amount of
generosity and charm! I would like to greet Volodya with his 60-th anniversary, 
and wish him health, new brilliant mathematical ideas, 
as well as many days filled with a positive spirit.

~\\[-0.3cm]

\noindent
\textbf{Acknowledgements:} I acknowledge NSF grants DMS-1161867 
and DMS-1501001 for the partial support. 
My work benefitted from participation in the program 
``Grothendieck-Teichm\"uller Groups, Deformation and Operads''
of the Isaac Newton Institute in Cambridge, UK. 
I would like to thank the organizers of this wonderful and 
very inspiring program. I would like thank  Rina Anno and Alexei Oblomkov
for asking me the question  ``Can we ask a computer to 
perform deformation quantization?''  The main theorem of this paper 
may be considered as the answer ``yes'' to this question.
The result of this paper was presented at the Geometry and Physics Seminar at
Boston University, the Deformation Theory seminar at Penn, and the Topology and Geometry Seminar 
at the Hebrew University of Jerusalem. 
I am thankful to the participants of these seminars for their questions and comments. 
A part of this text was written when I was a visitor at 
the Weizmann Institute of Science in Israel. I would like to thank the Weizmann Institute 
for hospitality and for wonderful conditions for visitors with families. I would like 
to thank Orit Dolgushev for her help with editing a couple of paragraphs of this paper. 
I would also like to thank anonymous referees for their comments about various 
versions of this paper.

\section{Reminder of $\OC$, $\KGra$ and stable formality quasi-isomorphisms}
\label{sec:OC-KGra-SFQ}

The definition of a stable formality quasi-isomorphism (SFQ) \cite[Section 5]{stable}  is 
based on the $2$-colored operads $\OC$ and $\KGra$\,. The 2-colored dg operad $\OC$ 
governs open-closed homotopy algebras introduced in \cite{OCHA} by  H. Kajiura and J. Stasheff
and the 2-colored operad $\KGra$ is ``assembled'' from graphs used in Kontsevich's paper \cite{K}.

\subsection{The Kajiura-Stasheff operad $\OC$}
As an operad in the category $\grVect$ of graded vector spaces, $\OC$ 
is freely generated by the 2-colored collection $\oc$ with the following 
spaces: 
\begin{equation}
\label{oc-mc}
\oc(n,0)^{\mc} = \bs^{3-2n} \bbK\,, \qquad  n \ge 2\,, 
\end{equation}
\begin{equation}
\label{oc-mo}
\oc(0,k)^{\mo} = \bs^{2-k}\, \sgn_k \otimes \bbK[S_k]\,, \qquad  k \ge 2\,, 
\end{equation}
\begin{equation}
\label{oc-mo1}
\oc(n,k)^{\mo} = \bs^{2-2n-k}\,  \sgn_k \otimes \bbK[S_k]\,, \qquad n \ge 1\,,
\qquad k \ge 0\,, 
\end{equation}
where $\sgn_k$ is the sign representation of $S_k$\,.
The remaining spaces of the collection $\oc$ are zero. 

Following \cite[Section 4]{stable}, 
we represent generators of $\OC$ in $\oc(n,0)^{\mc}$
by non-planar labeled corollas with $n$ solid incoming edges 
(see figure \ref{fig:bt-n-oc-mc}). We represent  generators of $\OC$ in 
$\oc(0,k)^{\mo}$ by planar labeled corollas with  $k$ dashed incoming edges
(see figure \ref{fig:bt-k-oc-mo}). 
Finally, we use labeled 2-colored corollas with a planar structure given only  
on the dashed edges to represent generators of $\OC$ in $\oc(n,k)^{\mo}$ 
(see figure \ref{fig:bt-nk-oc-mo}).
\begin{figure}[htp] 
\begin{minipage}[t]{0.3\linewidth}
\centering 
\begin{tikzpicture}[scale=0.5]
\tikzstyle{w}=[circle, draw, minimum size=3, inner sep=1]
\tikzstyle{vt}=[circle, draw, fill, minimum size=0, inner sep=1]
\node[vt] (l1) at (0, 2) {};
\draw (0,2.5) node[anchor=center] {{\small $1$}};
\node[vt] (l2) at (1, 2) {};
\draw (1,2.5) node[anchor=center] {{\small $2$}};
\draw (2,1.8) node[anchor=center] {{\small $\dots$}};
\node[vt] (ln) at (3, 2) {};
\draw (3,2.5) node[anchor=center] {{\small $n$}};
\node[w] (v) at (1.5, 1) {};
\node[vt] (r) at (1.5, 0) {};
\draw (r) edge (v);
\draw (v) edge (l1);
\draw (v) edge (l2);
\draw (v) edge (ln);
\end{tikzpicture}
\caption{The non-planar corolla $\st^{\mc}_n$
representing a generator of $\oc(n,0)^{\mc}$ } \label{fig:bt-n-oc-mc}
\end{minipage}
\hspace{0.2cm}
\begin{minipage}[t]{0.3\linewidth}
\centering 
\begin{tikzpicture}[scale=0.5]
\tikzstyle{vt}=[circle, draw, fill, minimum size=0, inner sep=1]
\tikzstyle{w}=[circle, draw, minimum size=3, inner sep=1]
\node[vt] (l1) at (0, 2) {};
\draw (0,2.5) node[anchor=center] {{\small $1$}};
\node[vt] (l2) at (1, 2) {};
\draw (1,2.5) node[anchor=center] {{\small $2$}};
\draw (2,1.8) node[anchor=center] {{\small $\dots$}};
\node[vt] (lk) at (3, 2) {};
\draw (3,2.5) node[anchor=center] {{\small $k$}};
\node[w] (v) at (1.5, 1) {};
\node[vt] (r) at (1.5, 0) {};
\draw [dashed] (r) edge (v);
\draw [dashed] (v) edge (l1);
\draw [dashed] (v) edge (l2);
\draw [dashed] (v) edge (lk);
\end{tikzpicture}
\caption{The 2-colored planar corola $\st^{\mo}_k$ representing 
a generator of $\oc(0,k)^{\mo}$ } \label{fig:bt-k-oc-mo}
\end{minipage}
\hspace{0.2cm}
\begin{minipage}[t]{0.35\linewidth}
\centering 
\begin{tikzpicture}[scale=0.5]
\tikzstyle{vt}=[circle, draw, fill, minimum size=0, inner sep=1]
\tikzstyle{w}=[circle, draw, minimum size=3, inner sep=1]
\node[vt] (l1) at (0, 2) {};
\draw (0,2.5) node[anchor=center] {{\small $1$}};
\draw (1.5,1.8) node[anchor=center] {{\small $\dots$}};
\node[vt] (ln) at (2, 2) {};
\draw (2,2.5) node[anchor=center] {{\small $n$}};
\node[vt] (l1mo) at (3, 2) {};
\draw (3,2.5) node[anchor=center] {{\small $1$}};
\draw (3.5,1.8) node[anchor=center] {{\small $\dots$}};
\node[vt] (lkmo) at (5, 2) {};
\draw (5,2.5) node[anchor=center] {{\small $k$}};
\node[w] (v) at (2.5, 1) {};
\node[vt] (r) at (2.5, 0) {};
\draw  [dashed] (r) edge (v);
\draw (v) edge (l1);
\draw (v) edge (ln);
\draw  [dashed] (v) edge (l1mo);
\draw  [dashed] (v) edge (lkmo);
\end{tikzpicture}
\caption{The 2-colored partially planar corola $\st^{\mo}_{n,k}$ representing 
a generator of $\oc(n,k)^{\mo}$ 
} \label{fig:bt-nk-oc-mo}
\end{minipage}
\end{figure}
Applying element $\si \in S_k$ to the labeled corolla $\st^{\mo}_k$
depicted in figure \ref{fig:bt-k-oc-mo} we get a basis for the vector 
space $\oc(0,k)^{\mo}$\,. Similarly, applying elements of 
$(\id , \si) \in S_n \times S_k$ to the labeled 
corolla depicted in figure \ref{fig:bt-nk-oc-mo} we get a basis 
for the vector space $\oc(n,k)^{\mo}$\,.

The corollas  $\st^{\mc}_n$, $\st^{\mo}_k$ and  $\st^{\mo}_{n,k}$
carry the following degrees: 
\begin{eqnarray}
\label{deg-mc-n}
| \st^{\mc}_n |  = 3 - 2n &  n \ge 2\,, \\[0.3cm]
\label{deg-mo-k}
| \st^{\mo}_k |  = 2 - k  &  k \ge 2\,, \\[0.3cm]
\label{deg-mo-nk}
| \st^{\mo}_{n,k} |  = 2 -2n - k &  ~~~~n \ge 1,~ k \ge 0\,. 
\end{eqnarray}

The differential $\cD$ on $\OC$ is defined by the following 
equations\footnote{For a nice pictorial definition of 
this differential on $\OC$ we refer the reader to \cite[Section 4.1]{stable}.}
\begin{equation}
\label{cD-st-mc}
\cD ( \st^{\mc}_n) : = - \sum_{p=2}^{n-1} \sum_{\tau \in \Sh_{p, n-p}} 
(\tau, \id) \big( \st^{\mc}_{n-p+1} \circ_{1, \mc} \st^{\mc}_{p}  \big)\,.
\end{equation}
\begin{equation}
\label{cD-st-mo-k}
\cD ( \st^{\mo}_k) : = - 
\sum_{p=0}^{k-2} \sum_{q= p+2}^{k} 
(-1)^{p + (k-q) (q-p)} \,
\st^{\mo}_{p + (k-q)+1}  \circ_{p+1, \mo} \st^{\mo}_{q-p}\,.
\end{equation}
\begin{equation}
\label{cD-st-mo-nk}
\cD ( \st^{\mo}_{n,k}) : =  (-1)^k 
\sum_{p=2}^{n} \sum_{\tau \in \Sh_{p, n-p}}
(\tau, \id) \big( \st^{\mo}_{n-p+1, k} \circ_{1, \mc}
\st^{\mc}_{p} \big)
\end{equation}
$$
- \sum_{ \substack{ 0 \le p \le q \le k  \\[0.1cm] p + (k-q) \ge 1 } }
(-1)^{p + (k-q)(q-p)}\,  \big( \st^{\mo}_{p + (k-q)+1} \circ_{p+1, \mo} \st^{\mo}_{n, q-p} \big)
$$
$$
 - 
\sum_{r=1}^{n-1}
 \sum_{\substack{ \si \in \Sh_{r, n-r} \\[0.1cm] 0 \le p \le q \le k }}
(-1)^{p + (k-q)(q-p)}\, (\si, \id) \big( \st^{\mo}_{r, \, p + (k-q)+1} \circ_{p+1, \mo} \st^{\mo}_{n-r, q-p} \big)\,. 
$$
$$
 - 
 ~ \sum_{ 0 \le p, ~ p+2 \le q \le k } ~
(-1)^{p + (k-q)(q-p)}\, \big( \st^{\mo}_{n,\,  p + (k-q)+1} \circ_{p+1, \mo} \st^{\mo}_{q-p} \big).
$$

Since the right hand sides of equations \eqref{cD-st-mc}, \eqref{cD-st-mo-k}, 
and \eqref{cD-st-mo-nk} are quadratic in generators, we conclude that the collection
$\bsi \oc$ is a pseudo-cooperad and
$$
\OC = \Cobar(\oc^{\vee}),
$$ 
where $\oc^{\vee}$ is the 2-colored coaugmented cooperad obtained from $\bsi \oc$
via ``adjoining the counit''. 

Finally, we recall that algebras over $\OC$ (a.k.a. open-closed homotopy algebras) 
are pairs of cochain complexes $(\cV, \cA)$ with 
the following data\footnote{Recall that introducing a $\La\Lie_{\infty}$-structure on $\cV$ is 
equivalent to introducing an $L_{\infty}$-structure on $\bsi \cV$\,.}:
 
\begin{itemize}

\item A $\La\Lie_{\infty}$-structure on $\cV$, 

\item an $A_{\infty}$-structure on $\cA$, and 

\item a $\La\Lie_{\infty}$-morphism from $\cV$ to the 
Hochschild cochain complex $\Cbu(\cA)$ of $\cA$\,.

\end{itemize} 

\subsection{The operad $\dGra$ and its $2$-colored version $\KGra$ }

To define the operad $\dGra$, we introduce a collection of 
auxiliary sets $\{\dgra_n \}_{n \ge 1}$\,.

\begin{figure}[htp] 
\centering 
\begin{tikzpicture}[scale=0.5, >=stealth']
\tikzstyle{w}=[circle, draw, minimum size=3, inner sep=1]
\tikzstyle{b}=[circle, draw, fill, minimum size=3, inner sep=1]
\node [b] (b1) at (0,0) {};
\draw (-0.4,0) node[anchor=center] {{\small $1$}};
\node [b] (b3) at (2,0) {};
\draw (2, 0.5) node[anchor=center] {{\small $3$}};
\node [b] (b2) at (6,0) {};
\draw (5.9, 0.6) node[anchor=center] {{\small $2$}};
\node [b] (b4) at (2,2) {};
\draw (2, 2.6) node[anchor=center] {{\small $4$}};
\node [b] (b5) at (5,2) {};
\draw (5, 2.6) node[anchor=center] {{\small $5$}};
\draw [->] (b3) edge (b1);
\draw (1,0.4) node[anchor=center] {{\small $i$}};
\draw [->] (b3) ..controls (3,0.5) and (5,0.5) .. (b2);
\draw (4,0.8) node[anchor=center] {{\small $ii$}}; 
\draw [<-] (b3) ..controls (3,-0.5) and (5,-0.5) .. (b2);
\draw (4,-0.8) node[anchor=center] {{\small $iii$}}; 
\end{tikzpicture}
\caption{Roman numerals indicate that  
$(3,1) < (3,2) < (2,3)$ } \label{fig:exam}
\end{figure}
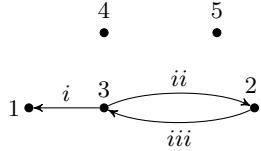

An element of $\dgra_{n}$ is a directed labelled graph $\G$
with $n$ vertices and with the additional piece 
of data: the set of edges of $\G$ is equipped with a 
total order. An example of an element in $\dgra_5$ is 
shown in figure \ref{fig:exam}. 
Here, we use roman numerals to specify a total order on a set of edges. 
For example, the roman numerals in figure \ref{fig:exam} indicate that  
$(3,1) < (3,2) < (2,3)$.

Next, we introduce a collection of graded vector spaces $\{\dGra(n) \}_{n \ge 1}$\,.
The space  $\dGra(n)$ is spanned by elements of 
$\dgra_n$, modulo the relation $\G^{\si} = (-1)^{|\si|} \G$\,,
where the graphs $\G^{\si}$ and $\G$ correspond to the same
directed labelled graph but differ only by permutation $\si$
of edges. We also declare that 
the degree of a graph $\G$ in $\dGra(n)$ equals 
$-e(\G)$, where $e(\G)$ is the number of edges in $\G$\,.
For example, the graph $\G$ in figure \ref{fig:exam} has 
$3$ edges. Thus its degree is $-3$\,. 

According to \cite{Thomas}, the collection $\{\dGra(n) \}_{n \ge 1}$
forms an operad in the category of graded vector spaces. 
The symmetric group $S_n$ acts on $\dGra(n)$ in the 
obvious way by rearranging labels and the operadic 
multiplications are defined in terms of natural operations 
of erasing vertices and attaching edges to vertices.

The operad $\dGra$ upgrades naturally to a 2-colored operad 
$\KGra$ whose spaces are finite linear combinations
of graphs used by M. Kontsevich in \cite{K}. 

For $\KGra$, we declare that $\KGra(n,k)^{\mc} = \bfzero  $ whenever $k \ge 1$\,. 

For the space $\KGra(n,0)^{\mc}$ ($n \ge 0$) we have 
\begin{equation}
\label{KGra-mc}
\KGra(n,0)^{\mc} = \dGra(n)\,.
\end{equation}

Finally, to define the space $\KGra(n,k)^{\mo}$ we introduce 
the auxiliary set $\dgra_{n,k}$\,.   
An element of the set $\dgra_{n,k}$ is a directed labelled graph $\G$
with $n$ vertices of color $\mc$, $k$ vertices of color $\mo$, and 
with the following data: the set of edges of $\G$ is equipped with a 
total order.
In addition, we require that  
each graph $\G \in \dgra_{n,k}$ has no edges originating from any vertex with 
color $\mo$\,. 

\begin{example}
\label{ex:KGra-mc-mo}
Figure \ref{fig:exam-mc-mo} shows an 
example of a graph in $\dgra_{2,3}$. Black (resp. white) vertices carry
the color $\mc$ (resp. $\mo$). We use separate labels for vertices of 
color $\mc$ and vertices of color $\mo$\,. For example $2_{\mc}$
denotes the vertex of color $\mc$ with label $2$ and 
$3_{\mo}$ denotes the vertex of color $\mo$ with label $3$\,.
\begin{figure}[htp] 
\centering 
\begin{tikzpicture}[scale=0.5, >=stealth']
\tikzstyle{w}=[circle, draw, minimum size=4, inner sep=1]
\tikzstyle{b}=[circle, draw, fill, minimum size=4, inner sep=1]
\node [w] (w2) at (0,0) {};
\draw (0,-0.6) node[anchor=center] {{\small $2$}};
\node [w] (w1) at (1,0) {};
\draw (1,-0.6) node[anchor=center] {{\small $1$}};
\node [w] (w3) at (2,0) {};
\draw (2,-0.6) node[anchor=center] {{\small $3$}};
\node [b] (b1) at (0,2) {};
\draw (0,2.6) node[anchor=center] {{\small $1$}};
\node [b] (b2) at (2,2) {};
\draw (2,2.6) node[anchor=center] {{\small $2$}};
\draw [->] (b1) edge (w1);
\draw [->] (b1) edge (b2);
\draw [->] (b2) edge (w1);
\draw [->] (b2) edge (w3);
\end{tikzpicture}
\caption{We equip the edges with the order 
$(1_{\mc}, 2_{\mc}) < (1_{\mc}, 1_{\mo}) <  (2_{\mc}, 1_{\mo})
<  (2_{\mc}, 3_{\mo}) $ } \label{fig:exam-mc-mo}
\end{figure}
\end{example}

The space $\KGra(n,k)^{\mo}$ is spanned by elements of 
$\dgra_{n,k}$, modulo the relation $\G^{\si} = (-1)^{|\si|} \G$\,, 
where the graphs $\G^{\si}$ and $\G$ correspond to the same
directed labelled graph
 but differ only by permutation $\si$
of edges. 
As above, we declare that 
the degree of a graph $\G$ in $\KGra(n,k)^{\mo}$ equals $-e(\G)$\,, 
where $e(\G)$ is the total number of edge of $\G$\,.

The operadic structure on the resulting $2$-colored collection $\KGra$
is defined in the similar way to that on $\dGra$. 
For more details, we refer the reader to \cite[Section 3]{stable}.

Just as in \cite{stable}, the following vectors of $\KGra$ will 
play a special role: 
%
%
\begin{equation}
\label{G-edge}
 \begin{tikzpicture}[scale=0.5, >=stealth']
\tikzstyle{w}=[circle, draw, minimum size=4, inner sep=1]
\tikzstyle{b}=[circle, draw, fill, minimum size=4, inner sep=1]
\draw (-2,0) node[anchor=center] {{$\G_{\ed} = $}};
\node [b] (b1) at (0,0) {};
\draw (0,0.6) node[anchor=center] {{\small $1$}};
\node [b] (b2) at (1.5,0) {};
\draw (1.5,0.6) node[anchor=center] {{\small $2$}};
\draw [->] (b1) edge (b2);
\draw (2.5,0) node[anchor=center] {{$+$}};
\node [b] (bb1) at (3.5,0) {};
\draw (3.5,0.6) node[anchor=center] {{\small $1$}};
\node [b] (bb2) at (5,0) {};
\draw (5,0.6) node[anchor=center] {{\small $2$}};
\draw [<-] (bb1) edge (bb2);
\end{tikzpicture}
\end{equation}
\begin{equation}
\label{binary}
\G_{\ww} =   \begin{tikzpicture}[scale=0.5, >=stealth']
\tikzstyle{w}=[circle, draw, minimum size=4, inner sep=1]
\tikzstyle{b}=[circle, draw, fill, minimum size=4, inner sep=1]
\node [w] (w1) at (0,0) {};
\draw (0,0.6) node[anchor=center] {{\small $1$}};
\node [w] (w2) at (1.5,0) {};
\draw (1.5,0.6) node[anchor=center] {{\small $2$}};
\draw (w1) (w2);
\end{tikzpicture}
\end{equation}

We also need the series of  ``brooms'' $\G^{\br}_k$ for $k \ge 0$ 
depicted in figure \ref{fig:brooms}.
\begin{figure}[htp]
\centering 
\begin{tikzpicture}[scale=0.5, >=stealth']
\tikzstyle{w}=[circle, draw, minimum size=4, inner sep=1]
\tikzstyle{b}=[circle, draw, fill, minimum size=4, inner sep=1]
\draw (-3,1.5) node[anchor=center] {{$\G^{\br}_{k}  \quad = \quad $}};
\node [b] (v) at (2.5,3) {};
\node [w] (v1) at (0,0) {};
\node [w] (v2) at (1,0) {};
\node [w] (vk) at (4,0) {};
\draw (2.5,0.7) node[anchor=center] {{\small $\dots$}};
\draw [->] (v) edge (1,1.2);
\draw  (1,1.2) edge (v1);
\draw [->] (v) edge (1.5,1);
\draw  (1.5,1) edge (v2);
\draw [->] (v) edge (3.5,1);
\draw  (3.5,1) edge (vk);
\draw (2.5,3.6) node[anchor=center] {{\small $1$}};
\draw (0,-0.6) node[anchor=center] {{\small $1$}};
\draw (1,-0.6) node[anchor=center] {{\small $2$}};
\draw (4,-0.6) node[anchor=center] {{\small $k$}};
\end{tikzpicture}
~\\[0.3cm]
\caption{Edges are ordered in this way 
$(1_{\mc}, 1_{\mo}) < (1_{\mc}, 2_{\mo}) < \dots < (1_{\mc}, k_{\mo})$
} \label{fig:brooms}
\end{figure} 

Note that the graph $\G^{\br}_{0} \in \KGra(1,0)^{\mo}$ 
consists of a single black vertex labeled by $1$ and it has no edges.  

%
%

\subsubsection{The projection $\Pi$ and the operator $\pa^{\Hoch}$}
\label{sec:Pi-and-Hoch}

Let us denote by 
\begin{equation}
\label{Pi-KGra}
\Pi\KGra(n,k)^{\mo}
\end{equation}
the subspace of all vectors in $\KGra(n,k)^{\mo}$ satisfying these properties: 
\begin{pty}
\label{P:one}
All white vertices in each graph of the 
linear combination $c$ have valency one.
\end{pty}
\begin{pty}
\label{P:anti-symm}
For every $\si \in S_k$ we have  
\begin{equation}
\label{app-anti-symm}
(\id , \si)\, (c) = (-1)^{|\si|} c\,.
\end{equation}
\end{pty} 
For example, the ``brooms'' $\G^{\br}_k$ depicted 
in figure \ref{fig:brooms} obviously satisfy these properties.
So $\G^{\br}_k \in \Pi\KGra(1,k)^{\mo}$\,.

For every vector $c \in  \KGra(n,k)^{\mo}$, we denote by $\Pi_1(c)$ 
the linear combination of graphs in $\dgra_{n,k}$ which is obtained from 
$c$ by retaining only graphs whose all white vertices are univalent.
The assignment 
\begin{equation}
\label{Pi-1}
c \mapsto \Pi_1(c) 
\end{equation}
is obviously a linear projection from $\KGra(n,k)^{\mo}$ onto the  
the subspace of vectors in $\KGra(n,k)^{\mo}$ which satisfy Property \ref{P:one}. 

Composing $\Pi_1$ with the alternation operator 
\begin{equation}
\label{Alt-mo}
\Alt^{\mo} =  \frac{1}{k!} \sum_{\si \in S_k} 
(-1)^{|\si|} (\id, \si)  ~ : ~  \KGra(n,k)^{\mo} \to \KGra(n,k)^{\mo}\,,
\end{equation}
we get a canonical projection  
\begin{equation}
\label{Pi}
\Pi := \Alt^{\mo} \circ \Pi_1 ~:~ \KGra(n,k)^{\mo} \to  \Pi\KGra(n,k)^{\mo}
\end{equation}
from $\KGra(n,k)^{\mo}$ onto $\Pi\KGra(n,k)^{\mo}$\,.
 
Just as in \cite{stable}, the following operator plays an important role:  
$$
\pa^{\Hoch} : \KGra(n,k)^{\mo} \to  \KGra(n,k+1)^{\mo}\,,
$$
%
%
\begin{equation}
\label{pa-Hoch}
\pa^{\Hoch} (\ga) = 
\G_{\ww} \,\c_{2, \mo}\, \ga - \ga \,\c_{1, \mo}\, \G_{\ww} + 
\ga \,\c_{2, \mo}\, \G_{\ww} - \dots 
\end{equation}
$$
+ (-1)^{k} \ga\, \c_{k, \mo}\, \G_{\ww}
+ (-1)^{k+1}\G_{\ww} \,\c_{1, \mo}\, \ga.
$$

It is not hard to see that $\pa^{\Hoch}$ commutes with the action 
of $S_n \times \{\id\}$, and $\pa^{\Hoch} \circ \pa^{\Hoch} = 0.$
So, using $\pa^{\Hoch}$, we can introduce the cochain complexes: 
\begin{equation}
\label{KGra-Hoch}
\KGra^{\Hoch} : = \bs^{2n-2} \, \bigoplus_{k \ge 0} \bs^{k} \KGra(n,k)^{\mo}\,, 
\quad 
 \KGra^{\Hoch}_{\inv} : =
\bs^{2n-2} \, \bigoplus_{k \ge 0} \bs^{k}\big( \KGra(n,k)^{\mo} \big)^{S_n}\,.
\end{equation} 

The cohomology groups of the complexes \eqref{KGra-Hoch} with the differential $\pa^{\Hoch}$ were
computed in \cite[Appendix A]{stable}. For example, due to \cite[Remark A.4]{stable}, every vector 
$$
c \in  \Pi\KGra(n,k)^{\mo}
$$
is $\pa^{\Hoch}$-closed. Moreover, 
$$
\Pi\KGra(n,k)^{\mo} \, \cap \, \pa^{\Hoch} \big( \KGra(n,k-1)^{\mo} \big) ~ = ~ \bfzero. 
$$

\subsection{Stable formality quasi-isomorphisms}
\label{sec:SFQ}

We recall from \cite{stable} that
\begin{defi} 
\label{dfn:stable}
A \emph{stable formality quasi-isomorphism} (SFQ) is a morphism of 2-colored dg operads
\begin{equation}
\label{eq:stable}
F :  \OC \to \KGra
\end{equation}
satisfying the following ``boundary conditions'': 
\begin{equation}
\label{Schouten}
F(\st^{\mc}_n) = \begin{cases}
   \G_{\ed} \qquad {\rm if} ~~ n = 2\,,  \\
   0 \qquad {\rm if} ~~  n \ge 3\,,
\end{cases}
\end{equation}
\begin{equation}
\label{mult}
F(\st^{\mo}_2) = \G_{\ww}\,, 
\end{equation}
and
\begin{equation}
\label{HKR}
F(\st^{\mo}_{1,k}) = \frac{1}{k!} \G^{\br}_{k}\,,
\end{equation}
where $\st^{\mc}_n$, $\st^{\mo}_k$, and $\st^{\mo}_{n,k}$ are corollas
depicted in figures \ref{fig:bt-n-oc-mc}, \ref{fig:bt-k-oc-mo}, \ref{fig:bt-nk-oc-mo}, 
respectively.
\end{defi}
 
Following \cite[Section 5.1]{stable}, we identify SFQs with MC elements $\al$
of the graded Lie algebra 
\begin{equation}
\label{Conv-oc-o-KGra}
\Conv(\bsi \oc, \KGra)
\end{equation}
subject to the three conditions: 
\begin{equation}
\label{al-Schouten}
\al(\bs^{-1}\, \st^{\mc}_n) = \begin{cases}
   \G_{\ed} \qquad {\rm if} ~~ n = 2\,,  \\
   0 \qquad {\rm if} ~~  n \ge 3\,,
\end{cases}
\end{equation}
\begin{equation}
\label{al-mult}
\al(\bs^{-1}\, \st^{\mo}_2) = \G_{\ww}\,, 
\end{equation}
and
\begin{equation}
\label{al-HKR}
\al(\bs^{-1} \, \st^{\mo}_{1,k}) = \frac{1}{k!} \G^{\br}_{k}\,.
\end{equation}

Let us observe that, since all vectors in $\KGra(0,k)^{\mo}$
have degree zero, we have 
\begin{equation}
\label{al-mult-higher}
\al(\bs^{-1}\, \st^{\mo}_k) = 0
\end{equation}
for all $k\ge 3$ and for all degree $1$ elements 
$\al$ in \eqref{Conv-oc-o-KGra}.

\begin{remark}
\label{rem:other-fields}
Although the operads $\OC$, $\KGra$ (as well as the 
graded Lie algebra \eqref{Conv-oc-o-KGra}) 
are defined over $\bbQ$, we often consider SFQs defined 
over a field extension $E$ of $\bbQ$. Such an SFQ is 
a morphism of 2-colored dg operads: 
$$
F :  \OC \otimes_{\bbQ} E \to \KGra  \otimes_{\bbQ} E
$$
satisfying \eqref{Schouten}, \eqref{mult}, \eqref{HKR}. 

Since $\OC$ is freely generated by the collection $\oc$ (as the operad in the category of 
graded vector spaces) and the vector spaces $\oc(n,k)^{\mc}$, $\oc(n,k)^{\mo}$
are finite dimensional, SFQs over $E$ are in bijection with MC elements $\al$ of 
the graded Lie algebra 
$$
\Conv(\bsi \oc, \KGra \otimes_{\bbQ} E )
$$
satisfying \eqref{al-Schouten}, \eqref{al-mult}, and \eqref{al-HKR}.

Since, in this paper, $E = \bbR$, we will often consider the graded Lie algebra 
\begin{equation}
\label{Conv-bbR}
\Conv(\bsi \oc, \KGra \otimes_{\bbQ} \bbR).
\end{equation} 
It goes without saying that every vector of \eqref{Conv-oc-o-KGra} may be viewed as
the vector of \eqref{Conv-bbR}. In particular, every SFQ defined over $\bbQ$ is 
an SFQ defined over $\bbR$ (i.e. a MC element of \eqref{Conv-bbR} satisfying 
\eqref{al-Schouten}, \eqref{al-mult}, and \eqref{al-HKR}).
\end{remark}
%
%
\begin{remark}
\label{rem:SFQ-filtration}
Let $\al$ be a degree $1$ element of \eqref{Conv-bbR} (or  \eqref{Conv-oc-o-KGra}) 
satisfying \eqref{al-Schouten}, \eqref{al-mult}, and \eqref{al-HKR}. It is not hard to see that 
\begin{equation}
\label{al-cF-mc-cF-mo}
\al \in \cF^{\mc}_{0} \Conv(\bsi \oc, \KGra \otimes_{\bbQ} \bbR), \qquad \al \in  
\cF^{\mo}_{-1}  \Conv(\bsi \oc, \KGra \otimes_{\bbQ} \bbR), 
\end{equation}
and 
$$
\al \notin \cF^{\mc}_{1} \Conv(\bsi \oc, \KGra \otimes_{\bbQ} \bbR), \qquad \al \notin
\cF^{\mo}_{0}  \Conv(\bsi \oc, \KGra \otimes_{\bbQ} \bbR), 
$$
where $\cF^{\mc}_{\bul}$ and $\cF^{\mo}_{\bul}$ are the descending filtrations on 
$\Conv(\bsi \oc, \KGra)$ corresponding to the colors $\mc$ and $\mo$
(see page \pageref{Conv-filtr-chi} for 
the definition of the filtration $\cF^{\chi}_{\bul}$).
\end{remark}

For our purposes, we will be interested in degree $1$ elements 
$\al \in \Conv(\bsi \oc, \KGra \otimes_{\bbQ} \bbR)$ (in particular, SFQs) 
satisfying the following technical property: 
%
%
\begin{pty}
\label{P:Pi-arity-2}
For every $k \ge 1$,
\begin{equation}
\label{Pi-alpha-2-k}
\Pi \big( \al(\bsi \st^{\mo}_{2, k})  \big) = 0,
\end{equation}
where $\Pi$ is defined in \eqref{Pi}.
\end{pty}

%
%

According to \cite{stable}, SFQs corresponding to MC elements
$\al$ and $\ti{\al}$ are homotopy equivalent if there exists a degree 
zero element 
$$
\xi \in \Conv(\bsi \oc, \KGra)
$$
(or more generally $\xi \in \Conv(\bsi \oc, \KGra \otimes_{\bbQ} \bbR)$) such that
\begin{equation}
\label{xi-cond-mc}
\xi(\bsi \st^{\mc}_n) =  0 \qquad \forall ~~ n \ge 2
\end{equation}
and
\begin{equation}
\label{homot-equiv}
\ti{\al} = \exp([\xi, ~])\, \al.
\end{equation}

\begin{remark}
\label{rem:xi-mo}
We would like to observe that, for every degree zero element 
$$
\xi  \in  \Conv(\bsi \oc, \KGra \otimes_{\bbQ} \bbR),
$$
\begin{equation}
\label{xi-mo}
\xi (\bsi \st^{\mo}_{k}) = 0, \qquad \forall~~ k \ge 2
\end{equation}
since the degree of $\bsi \st^{\mo}_{k}$ is $1-k < 0$ and every 
graph in $\dgra_{0,k}$ (i.e. without black vertices) has degree $0$.

Similarly, for every degree zero element $\xi \in  \Conv(\bsi \oc, \KGra \otimes_{\bbQ} \bbR)$, 
\begin{equation}
\label{xi-mo-1-k}
\xi (\bsi \st^{\mo}_{1, k}) = 0, \qquad \forall~~ k \ge 0.
\end{equation}
Indeed, the vector $\bsi \st^{\mo}_{1, k}$ has degree $-1-k$
and all graphs in $\dgra_{1,k}$ have $\le k$ edges. The latter is easy to 
see since multiple edges and loops are not allowed.  
\end{remark}
\begin{remark}
\label{rem:pty-Pi-stable}
Let $\beta$ be a MC element of \eqref{Conv-bbR} corresponding to an SFQ 
and satisfying Property \ref{P:Pi-arity-2}. 
Using the filtration $\cF^{\mc}_{\bul}$ and equation \eqref{xi-mo-1-k}, 
it is not hard to show that
$$
(e^{\ad_{\xi}}(\beta) - \beta) \, \big(\bsi \st^{\mo}_{2,k}) = - \pa^{\Hoch} \xi(\bsi \st^{\mo}_{2,k-1})
$$ 
for every degree $0$ vector $\xi$ in \eqref{Conv-bbR} satisfying \eqref{xi-cond-mc}.
Therefore identity $\Pi \circ \pa^{\Hoch} = 0$ implies that  
$$
\Pi \big( \, (e^{\ad_{\xi}}(\beta) - \beta) \, \big(\bsi \st^{\mo}_{2,k})  \big) = 0.  
$$
Thus Property \ref{P:Pi-arity-2} is stable under homotopy equivalences. 

In addition, it is not hard to see that Property \ref{P:Pi-arity-2} is stable under 
the action of the full directed graph complex $\dfGC$. (See \cite[Section 6]{stable}.)
\end{remark}

\subsection{Kontsevich's SFQ}
\label{sec:Konts}
The first example of a stable formality quasi-isomorphism 
(over $\bbR$) was constructed in \cite{K} by M. Kontsevich. 
In paper \cite{K}, M. Kontsevich did not use the language of operads. 
However this language is very convenient for our purposes.

Let us briefly recall here Kontsevich's construction of a particular example of an SFQ.
 
For a pair of integers $(n,k)$, $n \ge 0$, $k \ge 0$ satisfying the inequality 
$2n + k  \ge 2 $, we denote by $\Conf_{n,k}$ the configuration space of 
$n$ labeled points in the upper half plane and $k$ labeled points on the 
real line: 
\begin{multline}
\label{Conf-n-k}
\Conf_{n,k} : =  
\big\{  ( z_1, z_2, \dots, z_n; q_1, q_2, \dots, q_k) ~\big|~ z_i \in \bbC, ~~ \Im(z_i)  > 0, ~~
q_j \in \bbR, \\ 
z_{i_1} \neq z_{i_2}  ~ \textrm{ for } ~ i_1 \neq i_2\,, 
 ~~ q_{j_1} \neq q_{j_2}  ~ \textrm{ for } ~ j_1 \neq j_2 \big\}\,.  
\end{multline}

Let us denote by $G^{(1)}$ the 2-dimensional connected Lie group of 
the following transformations of the complex plane: 
\begin{equation}
\label{G-1}
G^{(1)} : = \{z \mapsto a z + b  ~|~ a, b \in \bbR, ~~ a > 0 \}\,.
\end{equation}

The condition $2n+k \ge 2$ guarantees that the diagonal action of 
$G^{(1)}$ on $\Conf_{n,k}$ is free and hence the quotient
\begin{equation}
\label{C-n-k}
C_{n,k} : =  \Conf_{n,k} ~\big/~ G^{(1)}
\end{equation}
is a smooth real manifold of dimension $2n+k - 2$\,.

We denote by $\bC_{n,k}$ the compactification of $C_{n,k}$ constructed 
by M. Kontsevich in \cite[Section 5]{K}. $\bC_{n,k}$ comes with an involved 
stratification which is described in great detail in {\it loc. cit.}  

Let $\G$ be a graph in $\dgra_{n,k}$ and $e$ be an edge of $\G$ which 
originates at the black\footnote{Let us recall that, by definition of $\dgra_{n,k}$ 
every edge $e$ of $\G \in \dgra_{n,k}$ should 
originate at a vertex with color $\mc$ (i.e. black vertex).} vertex with label $i$.
To such an edge $e$, we assign a 1-form $d \vf_{e}$
on $\Conf_{n,k}$ by the following rule: 
\begin{itemize}

\item if the edge $e$ terminates at the black vertex with label $j$ then 
$$
d \vf_e : =  d \Arg (z_j - z_i) - d \Arg  (z_j - \bar{z}_i),
$$

\item if the edge $e$ terminates at the white vertex with label $j$ then
$$
d \vf_e : =  d \Arg (q_j - z_i) - d \Arg  (q_j - \bar{z}_i) = 2 \, d \Arg (q_j - z_i). 
$$
 
\end{itemize}

It is easy to see that $d \vf_e$ descends to a 1-form on $C_{n,k}$. 
Furthermore, $d \vf_e$ extends to a smooth $1$-form on Kontsevich's 
compactification $\bC_{n,k}$ of $C_{n,k}$.

Using the embedding 
$$
\Conf_{n,k} \subset \bbC^n \times \bbR^k
$$
we equip the manifold $\Conf_{n,k}$ with the natural orientation 
which descends to $C_{n,k}$ and extends to $\bC_{n,k}$. 

To every element $\G \in \dgra_{n,k}$, we assign the following weight: 
\begin{equation}
\label{W-G}
W_{\G} : = \frac{1}{(2\pi)^{2n + k -2}} ~\int_{\bC^+_{n,k}}~ \bigwedge_{e \in E(\G)}  d \vf_{e}\,, 
\end{equation}
where $\bC^+_{n,k}$ is the closure of the connected component $C^+_{n,k}$ of
$C_{n,k}$ formed by configurations satisfying the condition
$$
q_1 < q_2 < \dots < q_k\,.
$$
The order of $1$-forms in $\displaystyle \bigwedge_{e \in E(\G)}  d \vf_{e}$ in \eqref{W-G}
agrees with the total order on the set of edges of $\G$\,. It is clear that the weight 
$W_{\G}$ for  $\G \in \dgra_{n,k}$ is non-zero only if the total number of 
edges of $\G$ equals $2n+ k - 2$.

Let us observe that the set $\dgra_{n,k}$ carries an obvious equivalence 
relation: two elements $\G$ and $\G'$ in $\dgra_{n,k}$ are equivalent if and 
only if they have the same underlying labeled (colored) graph. We denote by 
$$
[\dgra_{n,k}] 
$$
the set of corresponding equivalence classes.

Finally we define a degree $1$ element $\beta^K \in \Conv(\bsi \oc, \KGra)$ by 
setting
\begin{equation}
\label{al-K-Schouten}
\beta^K(\bs^{-1}\, \st^{\mc}_n) : = \begin{cases}
   \G_{\ed} \qquad {\rm if} ~~ n = 2\,,  \\
   0 \qquad {\rm if} ~~  n \ge 3\,,
\end{cases}
\end{equation}
\begin{equation}
\label{al-K-mult}
\beta^K(\bs^{-1}\, \st^{\mo}_2) : = \G_{\ww}\,, 
\end{equation}
and 
\begin{equation}
\label{al-K-main}
\beta^K (\bsi \st^{\mo}_{n,k}) : = \sum_{\ka  \in [\dgra_{n,k} ] } W_{\G_{\ka}} \, \G_{\ka}\,,
\end{equation}
where $\G_{\ka}$ is any representative of the equivalence class $\ka \in [\dgra_{n,k}]$\,.

Notice that, since the vector 
$$
W_{\G} \, \G ~\in~ \KGra(n,k)^{\mo}
$$
depends only on the equivalence class of $\G$\,, the right hand side 
of \eqref{al-K-main} does not depend on the choice of representatives 
$\G_{\ka}$\,.

A direct computation shows that the weights of the ``brooms''
depicted in figure \ref{fig:brooms} are given by the formula
$$
W_{\G^{\br}_k}  = \frac{1}{k!}\,.
$$
Hence $\beta^K$ satisfies all the required ``boundary conditions''
\eqref{al-Schouten}, \eqref{al-mult}, \eqref{al-HKR}\,. 

Following the line of arguments of \cite[Section 6.4]{K} one can 
show that $\beta^K$ satisfies the MC equation 
$$
[\beta^K, \beta^K] = 0 
$$
in \eqref{Conv-bbR}.
Thus $\beta^K$ gives us an SFQ. 
\begin{remark}
\label{rem:Weight}
The weight of a graph $\G \in \dgra_{n,k}$ defined in \cite[Section 6.2]{K}
comes with additional factors. These factors are absent in \eqref{W-G}
because our identification between polyvector fields and ``functions'' on the odd 
cotangent bundle is different from the one used by M. Kontsevich 
in \cite[Section 6.3]{K}.
\end{remark}

In Appendix \ref{app:Konts-ptys}, we prove that Kontsevich's SFQ $\beta^{K}$ satisfies 
Property \ref{P:Pi-arity-2}.

\section{Preparation}
\label{sec:prep}

Let $\al$ be a degree $1$ element of the graded Lie algebra \eqref{Conv-bbR}
satisfying \eqref{al-Schouten}, \eqref{al-mult} and \eqref{al-HKR}. Recall that
equation \eqref{al-mult-higher} holds for every degree $1$ element $\al$. 
Thus any degree $1$ element $\al$ of the Lie algebra \eqref{Conv-bbR}
(or the Lie algebra \eqref{Conv-oc-o-KGra})
satisfying \eqref{al-Schouten}, \eqref{al-mult} and \eqref{al-HKR}
is uniquely determined by the vectors
$$
\al(\bsi \st^{\mo}_{m,k}) \in \KGra(m,k)^{\mo} \otimes_{\bbQ} \bbR
$$ 
for $m \ge 2$ and $k \ge 0$. 

The goal of the following omnibus proposition is to describe 
which values of $\al$ may show up in the expressions
$$
[\al, \beta](\bsi \st^{\mo}_{n,k}) \qquad \textrm{and} \qquad [\al, \al](\bsi \st^{\mo}_{n,k}),
$$
where $\beta$ is a (possibly different) degree $1$ element of \eqref{Conv-bbR}.  
%
%
\begin{prop}
\label{prop:al-al}
Let $\al, \beta$ be degree $1$ elements of \eqref{Conv-bbR} 
satisfying \eqref{al-Schouten}, \eqref{al-mult} and \eqref{al-HKR}.

\begin{enumerate}

\item If $(m,k)$ is a pair of integers such that $m \ge 2$ and $k \ge 2$,
then the expression
\begin{equation}
\label{brack-al-beta}
[\al, \beta](\bsi \st^{\mo}_{m,k}) 
\end{equation}
\emph{does not} involve $\al(\bsi \st^{\mo}_{m', k'})$ if $(m', k')$ lies outside of 
the subset (see figure \ref{fig:what-contributes})
\begin{equation}
\label{the-rect-area}
\{ (m, k-1) \} ~ \sqcup~ \{(m', k')~|~ 1 \le m' \le m-1 \quad \textrm{and} \quad 0 \le k' \le k+1\}
\end{equation}
of $\bbZ_{\ge 1} \times \bbZ_{\ge 0}$.

\item  If $m \ge 3$, then the vectors $\al(\bsi \st^{\mo}_{m, k-1})$, $\al(\bsi \st^{\mo}_{m-1, k})$, 
and $\al(\bsi \st^{\mo}_{m-1, k+1})$
enter the expression
\begin{equation}
\label{al-al-m-k}
[\al, \al](\bsi \st^{\mo}_{m,k})
\end{equation}
linearly. Moreover, the vector $\al(\bsi \st^{\mo}_{2, k-1})$ enters the expression 
\begin{equation}
\label{al-al-2-k}
[\al, \al](\bsi \st^{\mo}_{2,k})
\end{equation}
linearly.

\item The expression
\begin{equation}
\label{brack-al-beta-k1}
[\al, \beta](\bsi \st^{\mo}_{m,1}) 
\end{equation}
\emph{does not} involve $\al(\bsi \st^{\mo}_{m', k'})$ if $(m', k')$ lies outside of the 
subset (see figure \ref{fig:smaller-area}): 
\begin{equation}
\label{smaller-area}
\{ (m',k') ~|~ 1 \le m' \le m-1 \quad \textrm{and} \quad 0 \le k' \le 2 \}
\end{equation}
of $\bbZ_{\ge 1} \times \bbZ_{\ge 0}$. 

\item  If $m \ge 3$, then the vectors $\al(\bsi \st^{\mo}_{m-1, 1})$ and $\al(\bsi \st^{\mo}_{m-1, 2})$
enter the expression
\begin{equation}
\label{al-al-m-1}
[\al, \al](\bsi \st^{\mo}_{m,1})
\end{equation}
linearly.

\item The expression
\begin{equation}
\label{brack-al-beta-k0}
[\al, \beta](\bsi \st^{\mo}_{m,0}) 
\end{equation}
\emph{does not} involve $\al(\bsi \st^{\mo}_{m', k'})$ if $(m', k')$ lies outside of the 
subset (see also figure \ref{fig:narrow-area}): 
\begin{equation}
\label{narrow-area}
\{(m',0) ~|~ 1 \le m' \le m-1\} ~\sqcup~ \{(m',1) ~|~ 1 \le m' \le m-1\}
\end{equation}
of $\bbZ_{\ge 1} \times \bbZ_{\ge 0}$. 

\item If $m \ge 3$, then the vectors $\al(\bsi \st^{\mo}_{m-1, 0})$ and $\al(\bsi \st^{\mo}_{m-1, 1})$
enter the expression
\begin{equation}
\label{al-al-m-0}
[\al, \al](\bsi \st^{\mo}_{m,0})
\end{equation}
linearly.

\item Finally, 
\begin{equation}
\label{al-beta-1-k}
[\al, \beta](\bsi \st^{\mo}_{1,k}) = 0, \qquad \forall ~~ k \ge 0. 
\end{equation}

\end{enumerate}

%
%
\begin{figure}[htp]
\centering 
\begin{tikzpicture}[scale=0.5, >=stealth']
\tikzstyle{bl}=[circle, draw, fill, minimum size=2, inner sep=1]
\tikzstyle{gr}=[circle, draw, fill, color=gray, minimum size=3, inner sep=1]
\draw [->] (-1,0) -- (7,0);
\draw [->] (0,-1) -- (0,10);
\draw (7, 0.4) node[anchor=center] {{\small $\mc$}};
\draw (0.4, 10) node[anchor=center] {{\small $\mo$}};
\draw (5, -0.5) node[anchor=center] {{\small $m$}};
\draw [thick] (5, -0.15) -- (5, 0.15);
\draw (1, -0.5) node[anchor=center] {{\small $1$}};
\draw (-0.5, 6) node[anchor=center] {{\small $k$}};
\draw [thick] (-0.15, 6) -- (0.15, 6);
\draw (-0.95, 7) node[anchor=center] {{\small $k+1$}};
\draw [thick] (-0.15, 7) -- (0.15, 7);
\draw [thick] (-0.1, 1) -- (0.1, 1);
\draw (-0.5, 1) node[anchor=center] {{\small $1$}};
\node [gr] at (5,6) {};
\node [bl] at (1,0) {};
\node [bl] at (2,0) {};
\node [bl] at (3,0) {};
\node [bl] at (4,0) {};
\node [bl] at (1,1) {};
\node [bl] at (2,1) {};
\node [bl] at (3,1) {};
\node [bl] at (4,1) {};
\node [bl] at (1,2) {};
\node [bl] at (2,2) {};
\node [bl] at (3,2) {};
\node [bl] at (4,2) {};
\node [bl] at (1,3) {};
\node [bl] at (2,3) {};
\node [bl] at (3,3) {};
\node [bl] at (4,3) {};
\node [bl] at (1,4) {};
\node [bl] at (2,4) {};
\node [bl] at (3,4) {};
\node [bl] at (4,4) {};
\node [bl] at (1,5) {};
\node [bl] at (2,5) {};
\node [bl] at (3,5) {};
\node [bl] at (4,5) {};
\node [bl] at (5,5) {}; \draw (5,5) circle (0.3);
\node [bl] at (1,6) {};
\node [bl] at (2,6) {};
\node [bl] at (3,6) {};
\node [bl] at (4,6) {}; \draw (4,6) circle (0.3);
\node [bl] at (1,7) {};
\node [bl] at (2,7) {};
\node [bl] at (3,7) {};
\node [bl] at (4,7) {}; \draw (4,7) circle (0.3);
\end{tikzpicture}
\caption{The gray bullet denotes the pair $(m,k)$ and the black bullets denote the pairs $(m',k')$ 
for which $\al(\bsi \st^{\mo}_{m', k'})$ may contribute to expression \eqref{brack-al-beta}
The circled bullets denote the pairs $(m',k')$ for which  $\al(\bsi \st^{\mo}_{m', k'})$ enter 
the expression \eqref{al-al-m-k} linearly if $m \ge 3$} \label{fig:what-contributes}
\end{figure}
%
%

%
%
\begin{figure}[htp]
\centering 
\begin{minipage}[t]{0.45\linewidth}
\centering 
\begin{tikzpicture}[scale=0.5, >=stealth']
\tikzstyle{bl}=[circle, draw, fill, minimum size=2, inner sep=1]
\tikzstyle{gr}=[circle, draw, fill, color=gray, minimum size=3, inner sep=1]
\draw [->] (-1,0) -- (7,0);
\draw [->] (0,-1) -- (0,4);
\draw (7,0.4) node[anchor=center] {{\small $\mc$}};
\draw (0.4,4) node[anchor=center] {{\small $\mo$}};
\draw (5, -0.5) node[anchor=center] {{\small $m$}};
\draw [thick] (5, -0.15) -- (5, 0.15);
\draw (1, -0.5) node[anchor=center] {{\small $1$}};
\draw (2, -0.5) node[anchor=center] {{\small $2$}};
\draw (-0.5, 2) node[anchor=center] {{\small $2$}};
\draw [thick] (-0.1, 2) -- (0.1, 2);
\draw (-0.5, 1) node[anchor=center] {{\small $1$}};
\draw [thick] (-0.1, 1) -- (0.1, 1);
\node [gr] at (5, 1) {};
\node [bl] at (1,0) {};
\node [bl] at (2,0) {};
\node [bl] at (3,0) {};
\node [bl] at (4,0) {};
\node [bl] at (1,1) {};
\node [bl] at (2,1) {};
\node [bl] at (3,1) {};
\node [bl] at (4,1) {}; \draw (4,1) circle (0.3);
\node [bl] at (1,2) {};
\node [bl] at (2,2) {};
\node [bl] at (3,2) {};
\node [bl] at (4,2) {};  \draw (4,2) circle (0.3);
\end{tikzpicture}
\caption{The gray bullet denotes the pair $(m,1)$ and the black bullets denote the pairs $(m',k')$ 
for which $\al(\bsi \st^{\mo}_{m', k'})$ 
may contribute to expression \eqref{brack-al-beta-k1}.
The two circled bullets denote the pairs $(m',k')$ for which $\al(\bsi \st^{\mo}_{m', k'})$ enter 
the expression \eqref{al-al-m-1} linearly if $m \ge 3$ } \label{fig:smaller-area}
\end{minipage}
\begin{minipage}[t]{0.05\linewidth}
~~~
\end{minipage}
\begin{minipage}[t]{0.45\linewidth}
\centering 
\begin{tikzpicture}[scale=0.5, >=stealth']
\tikzstyle{bl}=[circle, draw, fill, minimum size=2, inner sep=1]
\tikzstyle{gr}=[circle, draw, fill, color=gray, minimum size=3, inner sep=1]
\draw [->] (0,-1) -- (0,4);
\draw [->] (-1,0) -- (7,0);
\draw (7, 0.4) node[anchor=center] {{\small $\mc$}};
\draw (0.4, 4) node[anchor=center] {{\small $\mo$}};
\draw (5, -0.5) node[anchor=center] {{\small $m$}};

\draw (1, -0.5) node[anchor=center] {{\small $1$}};
\draw (2, -0.5) node[anchor=center] {{\small $2$}};
\draw (-0.5, 1) node[anchor=center] {{\small $1$}};
\draw [thick] (-0.1, 1) -- (0.1, 1);
\node [gr] at (5,0) {};
\node [bl] at (1,0) {};
\node [bl] at (2,0) {};
\node [bl] at (3,0) {};
\node [bl] at (4,0) {};  \draw (4,0) circle (0.3);
\node [bl] at (1,1) {};
\node [bl] at (2,1) {};
\node [bl] at (3,1) {};
\node [bl] at (4,1) {};  \draw (4,1) circle (0.3);
\end{tikzpicture}
\caption{The gray bullet denotes the pair $(m,0)$ and the black bullets denote the pairs $(m',k')$ 
for which $\al(\bsi \st^{\mo}_{m', k'})$ 
may contribute to expression \eqref{brack-al-beta-k0}.
The two circled bullets denote the pairs $(m',k')$ for which $\al(\bsi \st^{\mo}_{m', k'})$ 
enter the expression \eqref{al-al-m-0} linearly if $m \ge 3$ } \label{fig:narrow-area}
\end{minipage}
\end{figure}
\end{prop}
\begin{proof}
Let us observe that, due to Remark \ref{rem:SFQ-filtration}, 
$$
\beta \in \cF^{\mc}_{0} \Conv(\bsi \oc, \KGra \otimes_{\bbQ} \bbR), \qquad \beta \in  
\cF^{\mo}_{-1}  \Conv(\bsi \oc, \KGra \otimes_{\bbQ} \bbR). 
$$

Hence, for every $(m,k) \in \bbZ_{\ge 1} \times \bbZ_{\ge 0}$, the expression
\begin{equation}
\label{al-beta-m-k}
[\al, \beta](\bsi \st^{\mo}_{m,k}) 
\end{equation}
\emph{does not} involve $\al(\bsi \st^{\mo}_{m', k'})$ if $(m', k')$ lies outside of 
the subset
\begin{equation}
\label{bigger-one}
\{(m', k') \in  \bbZ_{\ge 1} \times \bbZ_{\ge 0} ~|~ m' \le m, \quad \textrm{and} \quad k' \le k+1 \}.
\end{equation}

According to \eqref{cD-st-mo-nk}, the expression 
$$
\cD (\st^{\mo}_{m,k}) 
$$
involves neither $\st^{\mo}_{m,k+1}$ nor $\st^{\mo}_{m,k}$. 
Moreover, since the elements $\st^{\mo}_{m,k'}$ with $k' < k$ show up only in the sums
$$
\sum_{ \substack{ 0 \le p \le q \le k  \\[0.1cm] p + (k-q) \ge 1 } }
(-1)^{p + (k-q)(q-p)}\,  \big( \st^{\mo}_{p + (k-q)+1} \circ_{p+1, \mo} \st^{\mo}_{m, q-p} \big),
$$
$$
\sum_{ 0 \le p, ~ p+2 \le q \le k } ~
(-1)^{p + (k-q)(q-p)}\, \big( \st^{\mo}_{m,\,  p + (k-q)+1} \circ_{p+1, \mo} \st^{\mo}_{q-p} \big),
$$
and (see \eqref{al-mult-higher}) 
$$
\beta(\bs^{-1}\, \st^{\mo}_k) = 0, \qquad \forall~~ k \ge 3,
$$
we see that, for $k' < k-1$, $\al(\bsi \st^{\mo}_{m,k'})$ does not show up 
in \eqref{al-beta-m-k}. 

Thus, in order to settle statements 1, 3, and 5, it remain to show that 
$\al(\bsi \st^{\mo}_{m,0})$ does not contribute to 
$$
[\al, \beta](\bsi \st^{\mo}_{m,1}).  
$$

Collecting the terms with $\st^{\mo}_{m,0}$ in $\cD( \st^{\mo}_{m,1})$, we get 
$$
\st^{\mo}_{2} \circ_{2,\mo} \st^{\mo}_{m,0} - \st^{\mo}_{2} \circ_{1,\mo} \st^{\mo}_{m,0}\,.
$$
Therefore, since 
$$
\G_{\ww}  \circ_{2, \mo} v -  \G_{\ww}  \circ_{1, \mo} v = 0, \qquad \forall ~~ v \in \KGra(m,0)^{\mo}\,,
$$
we see that $\al(\bsi \st^{\mo}_{m,0})$ indeed does not contribute to $[\al, \beta](\bsi \st^{\mo}_{m,1})$.  

Thus statements 1, 3, and 5 of the proposition are proved. 

Let us now assume that $m \ge 3$ and consider the expression \eqref{al-al-m-k}. 
It is obvious that  $\al(\bsi \st^{\mo}_{m, k-1})$ enter this expression linearly. 

Moreover, the terms in $\cD(\st^{\mo}_{m,k})$ involving the corollas $\st^{\mo}_{m-1, k+1}$ 
and $\st^{\mo}_{m-1, k}$ are of the form 
$$
(\tau, \id) \big( \st^{\mo}_{m-1, k} \circ_{1, \mc} \st^{\mc}_2 \big)\,,  \qquad
(\tau, \id) \big( \st^{\mo}_{m-1, k+1} \circ_{i, \mo} \st^{\mo}_{1,0} \big)\,,  
$$
$$
(\tau, \id) \big( \st^{\mo}_{m-1, k} \circ_{j, \mo} \st^{\mo}_{1,1} \big)\,, \qquad 
(\tau, \id) \big( \st^{\mo}_{1,1} \circ_{1, \mo} \st^{\mo}_{m-1, k} \big)\,, 
$$
where $\tau$ is a permutation in $S_m$, $ 1 \le i \le k+1$ and $1 \le j \le k$.  

Thus, since $m-1 \ge 2$, the vectors $\al(\bsi \st^{\mo}_{m-1, k+1})$ and 
$\al(\bsi \st^{\mo}_{m-1, k})$ indeed enter the expression \eqref{al-al-m-k}
linearly. 

The expressions \eqref{al-al-2-k}, \eqref{al-al-m-1} and \eqref{al-al-m-0} are treated analogously. 

Let us finally prove statement 7 of the proposition. 

The cases $k=0$ and $k=1$ are straightforward. So we only consider $k \ge 2$.  

Since every degree $1$ vector in \eqref{Conv-bbR} satisfies \eqref{al-mult-higher}, it is 
easy to see that \eqref{al-beta-1-k} is equivalent to 
\begin{equation}
\label{Hoch-beta-al}
\pa^{\Hoch}\, \beta(\bsi \st^{\mo}_{1, k-1}) + \pa^{\Hoch}\, \al(\bsi \st^{\mo}_{1, k-1})  =0,
\end{equation}
where $\pa^{\Hoch}$ is defined in \eqref{pa-Hoch}.

Equation \eqref{Hoch-beta-al} is satisfied since 
$$
\beta(\bsi \st^{\mo}_{1, k'}), ~~ \al(\bsi \st^{\mo}_{1, k'})  \in \Pi\KGra(1,k')^{\mo}
$$
and every vector in $\Pi\KGra(1,k')^{\mo}$ is $\pa^{\Hoch}$-closed. 

Proposition \ref{prop:al-al} is proved. 
\end{proof}

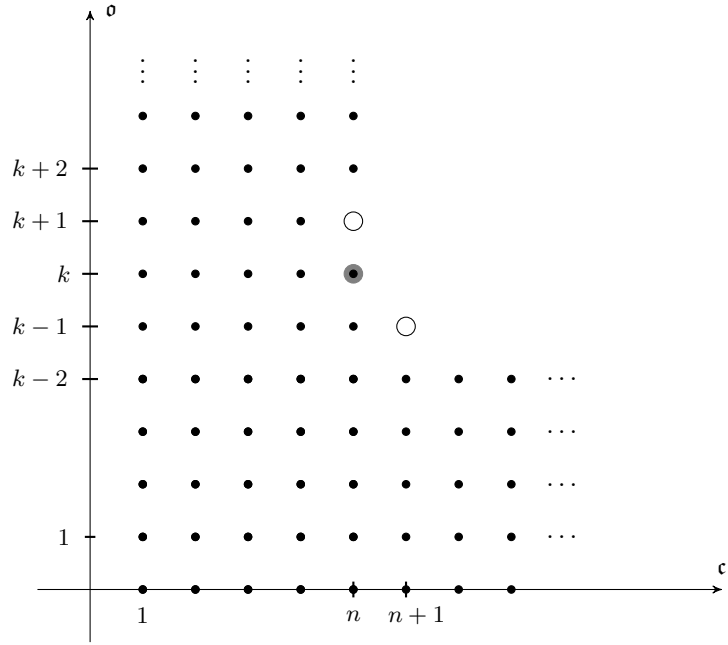
\begin{figure}[htp]
\centering 
\begin{tikzpicture}[scale=0.7, >=stealth']
\tikzstyle{bl}=[circle, draw, fill, minimum size=2, inner sep=1]
\tikzstyle{gr}=[circle, draw, fill, color=gray, minimum size=7, inner sep=1]
\tikzstyle{wh}=[circle, draw, minimum size=7, inner sep=1]
\draw [->] (-1,0) -- (12,0);
\draw [->] (0,-1) -- (0,11);
\draw (12, 0.4) node[anchor=center] {{\small $\mc$}};
\draw (0.4, 11) node[anchor=center] {{\small $\mo$}};
\draw (5, -0.5) node[anchor=center] {{\small $n$}};
\draw [thick] (5, -0.15) -- (5, 0.15);
\draw (6.2, -0.5) node[anchor=center] {{\small $n+1$}};
\draw [thick] (6, -0.15) -- (6, 0.15);
\draw (1, -0.5) node[anchor=center] {{\small $1$}};
\draw (-0.5, 6) node[anchor=center] {{\small $k$}};
\draw [thick] (-0.15, 6) -- (0.15, 6);
\draw (-0.95, 7) node[anchor=center] {{\small $k+1$}};
\draw [thick] (-0.15, 7) -- (0.15, 7);
\draw (-0.95, 8) node[anchor=center] {{\small $k+2$}};
\draw [thick] (-0.15, 8) -- (0.15, 8);
\draw (-0.95, 5) node[anchor=center] {{\small $k-1$}};
\draw [thick] (-0.15, 5) -- (0.15, 5);
\draw (-0.95, 4) node[anchor=center] {{\small $k-2$}};
\draw [thick] (-0.15, 4) -- (0.15, 4);
\draw [thick] (-0.1, 1) -- (0.1, 1);
\draw (-0.5, 1) node[anchor=center] {{\small $1$}};
\node [wh] at (5,7) {};
\node [wh] at (6,5) {};
\node [gr] at (5,6) {};
\node [bl] at (5,6) {};
\node [bl] at (1,0) {};
\node [bl] at (2,0) {};
\node [bl] at (3,0) {};
\node [bl] at (4,0) {};
\node [bl] at (5,0) {};
\node [bl] at (6,0) {};
\node [bl] at (7,0) {};
\node [bl] at (8,0) {};
\node [bl] at (1,1) {};
\node [bl] at (2,1) {};
\node [bl] at (3,1) {};
\node [bl] at (4,1) {};
\node [bl] at (5,1) {};
\node [bl] at (6,1) {};
\node [bl] at (7,1) {};
\node [bl] at (8,1) {};
\node [bl] at (1,2) {};
\node [bl] at (2,2) {};
\node [bl] at (3,2) {};
\node [bl] at (4,2) {};
\node [bl] at (5,2) {};
\node [bl] at (6,2) {};
\node [bl] at (7,2) {};
\node [bl] at (8,2) {};
\node [bl] at (1,3) {};
\node [bl] at (2,3) {};
\node [bl] at (3,3) {};
\node [bl] at (4,3) {};
\node [bl] at (5,3) {};
\node [bl] at (6,3) {};
\node [bl] at (7,3) {};
\node [bl] at (8,3) {};
\node [bl] at (1,4) {};
\node [bl] at (2,4) {};
\node [bl] at (3,4) {};
\node [bl] at (4,4) {};
\node [bl] at (5,4) {};
\node [bl] at (6,4) {};
\node [bl] at (7,4) {};
\node [bl] at (8,4) {};
\node [bl] at (1,0) {};
\node [bl] at (1,1) {};
\node [bl] at (1,2) {};
\node [bl] at (1,3) {};
\node [bl] at (1,4) {};
\node [bl] at (1,5) {};
\node [bl] at (1,6) {};
\node [bl] at (1,7) {};
\node [bl] at (1,8) {};
\node [bl] at (1,9) {};
\node [bl] at (2,0) {};
\node [bl] at (2,1) {};
\node [bl] at (2,2) {};
\node [bl] at (2,3) {};
\node [bl] at (2,4) {};
\node [bl] at (2,5) {};
\node [bl] at (2,6) {};
\node [bl] at (2,7) {};
\node [bl] at (2,8) {};
\node [bl] at (2,9) {};
\node [bl] at (3,0) {};
\node [bl] at (3,1) {};
\node [bl] at (3,2) {};
\node [bl] at (3,3) {};
\node [bl] at (3,4) {};
\node [bl] at (3,5) {};
\node [bl] at (3,6) {};
\node [bl] at (3,7) {};
\node [bl] at (3,8) {};
\node [bl] at (3,9) {};
\node [bl] at (4,0) {};
\node [bl] at (4,1) {};
\node [bl] at (4,2) {};
\node [bl] at (4,3) {};
\node [bl] at (4,4) {};
\node [bl] at (4,5) {};
\node [bl] at (4,6) {};
\node [bl] at (4,7) {};
\node [bl] at (4,8) {};
\node [bl] at (4,9) {};
\node [bl] at (5,0) {};
\node [bl] at (5,1) {};
\node [bl] at (5,2) {};
\node [bl] at (5,3) {};
\node [bl] at (5,4) {};
\node [bl] at (5,5) {};
\node [bl] at (5,8) {}; 
\node [bl] at (5,9) {}; 
%
\draw (9, 4) node[anchor=center] {{$\dots$}};
\draw (9, 3) node[anchor=center] {{$\dots$}};
\draw (9, 2) node[anchor=center] {{$\dots$}};
\draw (9, 1) node[anchor=center] {{$\dots$}};
\draw (1, 10) node[anchor=center] {{$\vdots$}};
\draw (2, 10) node[anchor=center] {{$\vdots$}};
\draw (3, 10) node[anchor=center] {{$\vdots$}};
\draw (4, 10) node[anchor=center] {{$\vdots$}};
\draw (5, 10) node[anchor=center] {{$\vdots$}};
\end{tikzpicture}
\caption{The big gray bullet is the point $(n,k)$. Points of the subset \eqref{area-no-change} are 
indicated by black bullets. The white circles indicate the points related to 
\eqref{n-k1-Hoch} and \eqref{n1-1k-broom0}} \label{fig:no-change}
\end{figure}

\begin{prop}
\label{prop:exp-xi-beta}
Let $(n,k)$ be a point in $\bbZ_{\ge 2} \times \bbZ_{\ge 1}$ and
$\beta$ be a MC element of  \eqref{Conv-bbR}
satisfying \eqref{al-Schouten}, \eqref{al-mult} and \eqref{al-HKR}. 
In addition, let $\xi$ be a degree zero element of \eqref{Conv-bbR} 
which satisfies \eqref{xi-cond-mc} and such that
\begin{equation}
\label{xi-only-for-n-k}
\xi(\bsi \st^{\mo}_{n_1,k_1}) = 0, \qquad \forall ~~ (n_1, k_1) \neq (n,k).
\end{equation}
Then 
\begin{equation}
\label{no-change}
\exp(\ad_{\xi})(\beta)\, (\bsi \st^{\mo}_{n_1, k_1}) = \beta (\bsi \st^{\mo}_{n_1, k_1})  
\end{equation}
for all pairs $(n_1, k_1)$ in the union\footnote{See figure \ref{fig:no-change} for the graphical 
presentation of this subset in  $\bbZ_{\ge 1} \times \bbZ_{\ge 0}$.}
\begin{equation}
\label{area-no-change}
\{(n_1, k_1) ~|~ n_1 < n \, \wedge \, k_1 \ge 0 \}
~\cup~
\{(n, k_1) ~|~  k_1 \neq k + 1  \} 
\end{equation}
$$
~\cup~ \{(n_1, k_1) ~|~ k_1 \le k-2  \, \wedge \, n_1 \ge 1 \}. 
$$
For the points $(n, k+1)$ and $(n+1, k-1)$
(indicated by white circles in figure \ref{fig:no-change}), 
we have
\begin{equation}
\label{n-k1-Hoch}
\exp(\ad_{\xi})(\beta)\, (\bsi \st^{\mo}_{n, k+1}) ~ = ~ \beta (\bsi \st^{\mo}_{n, k+1})  - 
\pa^{\Hoch} \xi (\bsi \st^{\mo}_{n, k}). 
\end{equation}

$$
\exp(\ad_{\xi})(\beta)\, (\bsi \st^{\mo}_{n+1, k-1}) ~ = ~  \beta (\bsi \st^{\mo}_{n+1, k-1}) 
$$
\begin{equation}
\label{n1-1k-broom0}
 ~ + ~
\sum_{i=1}^{n+1} \sum_{0 \le p \le k-1} 
(-1)^p \big(\tau_{n+1, i}, \id \big) \big(\, \xi (\bsi \st^{\mo}_{n, k}) \circ_{p+1, \mo} \G^{\br}_0 \, \big),
\end{equation}
where $\tau_{n+1,i}$ is the cycle $(i, i+1, \dots, n, n+1)$ in $S_{n+1}$. 
\end{prop}
%
%
\begin{proof} Let us use the filtration $\cF^{\mc}_{\bul}$ and $\cF^{\mo}_{\bul}$ 
on $\Conv(\bsi \oc, \KGra \otimes_{\bbQ} \bbR)$ (see page \pageref{Conv-filtr-chi} for 
the definition of the filtration $\cF^{\chi}_{\bul}$). 

Due to the first inclusion in \eqref{al-cF-mc-cF-mo} and 
$$
\xi \in \cF^{\mc}_{n} \Conv(\bsi \oc, \KGra \otimes_{\bbQ} \bbR),
$$
we have
$$
\ad^q_{\xi} (\beta)  \in \cF^{\mc}_{n} \Conv(\bsi \oc, \KGra \otimes_{\bbQ} \bbR)
$$
for all $q \ge 1$. Hence \eqref{no-change} holds if $n_1 < n$. 

Let us assume that $k \ge 2$. 
Due to the second inclusion in \eqref{al-cF-mc-cF-mo} and  
$$
\xi \in \cF^{\mo}_{k-1} \Conv(\bsi \oc, \KGra \otimes_{\bbQ} \bbR),
$$
we have  
$$
\ad^q_{\xi} (\beta)  \in \cF^{\mo}_{k-2} \Conv(\bsi \oc, \KGra \otimes_{\bbQ} \bbR)
$$
for all $q \ge 1$. Hence \eqref{no-change} holds if $k_1 \le k - 2$. 

Let us now consider $\ad^q_{\xi} (\beta) (\bsi \st^{\mo}_{n, k_1}) $ for $q \ge 1$. 
If $k_1 \le k$ then all tensor factors in 
$$
\cD(\st^{\mo}_{n, k_1} )
$$
of the form $\st^{\mo}_{n_2, k_2}$ have $n_2 < n$ or $k_2 < k$. 
Thus, using \eqref{xi-cond-mc} and \eqref{xi-mo}, we get that  
$$
\ad^q_{\xi} (\beta) (\bsi \st^{\mo}_{n, k_1}) =  0, \qquad \forall ~~ k_1 \le k, ~~ q \ge 1.
$$

To take care of $\ad^q_{\xi} (\beta) (\bsi \st^{\mo}_{n, k_1}) $ for $k_1 \ge k+2$, 
we observe that all terms in 
$$
\cD(\st^{\mo}_{n, k_1})
$$
which involve $\st^{\mo}_{n, k}$ are of the form 
$$
\st^{\mo}_{n, k} \circ_{p, \mo} \st^{\mo}_{\ti{k}}\,, \qquad 
\textrm{or} \qquad  \st^{\mo}_{\ti{k}} \circ_{p, \mo} \st^{\mo}_{n, k} 
$$
for $\ti{k} \ge 3$ and some $p$. Therefore, since
$$
\ad^{q-1}_{\xi} (\beta) (\bsi  \st^{\mo}_{\ti{k}})  = 0, \qquad \forall ~ \ti{k} \ge 3
$$
for the degree reason, we have 
$$
\ad^q_{\xi} (\beta) (\bsi \st^{\mo}_{n, k_1}) =  0, \qquad \forall ~~ k_1 \ge  k + 2, ~~ q \ge 1.
$$
Now it remains to prove \eqref{n-k1-Hoch} and \eqref{n1-1k-broom0}. 

Applying \eqref{brack-pa-cobar}, we get 
$$
[\xi, \beta] (\bsi \st^{\mo}_{n, k+1}) = - \pa^{\Hoch} \big( \xi( \bsi \st^{\mo}_{n,k} ) \big).
$$

Since only the sum 
$$
\st^{\mo}_{2} \circ_{2, \mo} \st^{\mo}_{n,k} - 
\sum_{p=0}^{k-1} (-1)^p  \st^{\mo}_{n,k} \circ_{p+1, \mo} \st^{\mo}_2
- (-1)^k \st^{\mo}_{2} \circ_{1, \mo} \st^{\mo}_{n,k}
$$
from $\cD(\st^{\mo}_{n, k+1})$ contributes to the bracket 
$$
[\xi, \ad^q_{\xi} (\beta)]  (\bsi \st^{\mo}_{n, k+1}), 
$$
and $\ad^q_{\xi} (\beta)\, (\bsi \st^{\mo}_2) = 0$ for all $q \ge 1$, we conclude that 
$$
[\xi, \ad^q_{\xi} (\beta)]  (\bsi \st^{\mo}_{n, k+1}) = 0, \qquad \forall ~~ q \ge 1.
$$
Thus \eqref{n-k1-Hoch} indeed holds. 

Finally, the proof of \eqref{n1-1k-broom0} is based on the direct 
application of \eqref{brack-pa-cobar}, 
\eqref{xi-cond-mc}, \eqref{xi-mo}, and \eqref{xi-only-for-n-k}.
\end{proof}

\section{The main theorem}
\label{sec:main-thm}

\subsection{Approximations to an SFQ} 
\label{sec:m-approx}
\begin{figure}[htp] 
\centering 
\begin{tikzpicture}[scale=1, >=stealth']
\tikzstyle{bl}=[circle, draw, fill, minimum size=4, inner sep=1]
\tikzstyle{gr}=[circle, draw, fill, color=gray, minimum size=7, inner sep=1]
\draw [->] (-1,0) -- (8,0);
\draw (8, -0.3) node[anchor=center] {{$n$}};
\draw [->] (0,-1) -- (0,13);
\draw (-0.3, 13) node[anchor=center] {{$k$}};
\draw[step=1, color=gray] (1,0) grid (7, 12); 
\draw [thick] (1, -0.1) -- (1,0.1);
\draw (1,-0.3) node[anchor=center] {{\small $1$}};
\draw [thick] (2, -0.1) -- (2,0.1);
\draw (2,-0.3) node[anchor=center] {{\small $2$}};
\draw [thick] (3, -0.1) -- (3,0.1);
\draw (3,-0.3) node[anchor=center] {{\small $3$}};
\draw [thick] (4, -0.1) -- (4,0.1);
\draw (4,-0.3) node[anchor=center] {{\small $4$}};
\draw [thick] (5, -0.1) -- (5,0.1);
\draw (5,-0.3) node[anchor=center] {{\small $5$}};
\draw [thick] (6, -0.1) -- (6,0.1);
\draw (6,-0.3) node[anchor=center] {{\small $6$}};
\draw [thick] (7, -0.1) -- (7,0.1);
\draw (7,-0.3) node[anchor=center] {{\small $7$}};
%
\draw [thick] (-0.1,1) -- (0.1,1);
\draw (-0.3,1) node[anchor=center] {{\small $1$}};
\draw [thick] (-0.1,2) -- (0.1,2);
\draw (-0.3,2) node[anchor=center] {{\small $2$}};
\draw [thick] (-0.1, 3) -- (0.1, 3);
\draw (-0.3, 3) node[anchor=center] {{\small $3$}};
\draw [thick] (-0.1, 4) -- (0.1, 4);
\draw (-0.3, 4) node[anchor=center] {{\small $4$}};
\draw [thick] (-0.1, 5) -- (0.1, 5);
\draw (-0.3, 5) node[anchor=center] {{\small $5$}};
\draw [thick] (-0.1, 6) -- (0.1, 6);
\draw (-0.3, 6) node[anchor=center] {{\small $6$}};
\draw [thick] (-0.1, 7) -- (0.1, 7);
\draw (-0.3, 7) node[anchor=center] {{\small $7$}};
\draw [thick] (-0.1, 8) -- (0.1, 8);
\draw (-0.3, 8) node[anchor=center] {{\small $8$}};
\draw [thick] (-0.1, 9) -- (0.1, 9);
\draw (-0.3, 9) node[anchor=center] {{\small $9$}};
\draw [thick] (-0.1, 10) -- (0.1, 10);
\draw (-0.3, 10) node[anchor=center] {{\small $10$}};
\draw [thick] (-0.1, 11) -- (0.1,11);
\draw (-0.3, 11) node[anchor=center] {{\small $11$}};
\draw [thick] (-0.1, 12) -- (0.1, 12);
\draw (-0.3, 12) node[anchor=center] {{\small $12$}};
%
%
\node [bl] at (1,0) {};
\node [bl] at (1,1) {}; 
\node [bl] at (1,2) {}; 
\node [bl] at (1,3) {}; 
\node [bl] at (1,4) {}; 
\node [bl] at (1,5) {}; 
\node [bl] at (1,6) {}; 
\node [bl] at (1,7) {}; 
\node [bl] at (1,8) {}; 
\node [bl] at (1,9) {}; 
\node [bl] at (1,10) {}; 
\node [bl] at (1,11) {}; 
\node [bl] at (1,12) {}; 
\node [bl] at (2,0) {};
\node [bl] at (2,1) {}; 
\node [bl] at (2,2) {}; 
\node [bl] at (2,3) {}; 
\node [bl] at (2,4) {}; 
\node [bl] at (2,5) {}; 
\node [bl] at (2,6) {}; 
\node [bl] at (2,7) {}; 
\node [bl] at (3,0) {};
\node [bl] at (3,1) {}; 
\node [bl] at (3,2) {}; 
\node [bl] at (3,3) {}; 
\node [bl] at (3,4) {}; 
\node [bl] at (3,5) {}; 
\node [bl] at (4,0) {};
\node [bl] at (4, 1) {}; 
\node [bl] at (4, 2) {}; 
\node [bl] at (4, 3) {}; 
\node [bl] at (5, 0) {};
\node [bl] at (5, 1) {}; 
%
%
\node [gr] at (6,0) {};
\node [gr] at (5,2) {};
\node [gr] at (4,4) {};
\node [gr] at (3,6) {};
\node [gr] at (2,8) {};
\end{tikzpicture}
\caption{In this example, $m=5$. $\A_m$ consists of the points shown 
as black bullets and $\B_m$ consists of the points shown as gray bullets} \label{fig:areas-m-approx}
\end{figure}
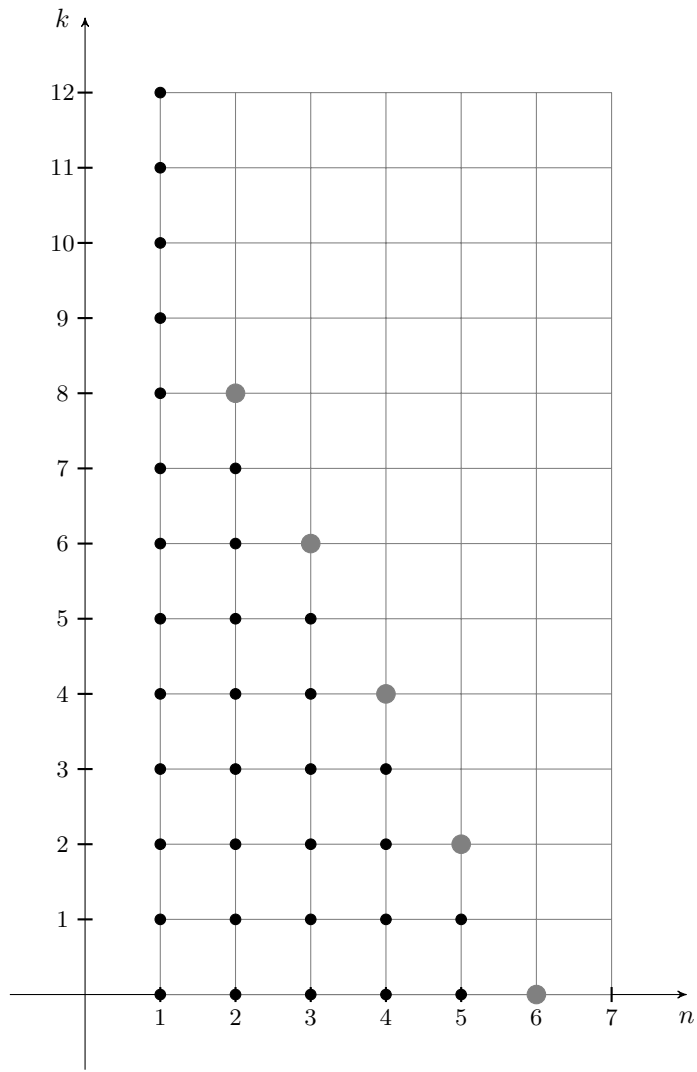

Let $m$ be an integer $\ge 2$,  $\A_{m}$ be the following 
subset of $\bbZ_{\ge 1} \times \bbZ_{\ge 0}$: 
\begin{equation}
\label{A-m-dfn}
\big\{ (n,k) \in \bbZ_{\ge 1} \times \bbZ_{\ge 0} ~ \big|~ 2n + k \le 2 m + 1   \big\} \cup
\big\{ (1,k)~ \big|~ \textrm{for any } k \ge 0  \big\},
\end{equation}
and $\B_m$ be this subset\footnote{Examples of $\A_m$ and $\B_m$,
for $m=5$, are depicted in figure \ref{fig:areas-m-approx}.} in $\bbZ_{\ge 1} \times \bbZ_{\ge 0}$
\begin{equation}
\label{B-m-dfn}
\B_m : = \big\{ (n, 2 (m-n) + 2) ~|~  2 \le n \le m+1 \big\},
\end{equation}
or equivalently
$$
\B_m : = \{(m+1, 0),\, (m, 2),\,  (m-1, 4),\, \dots,\, (2, 2m-2) \}.
$$

Let $\al$ be a degree $1$ element of $\Conv(\bsi \oc, \KGra)$
satisfying \eqref{al-Schouten}, \eqref{al-mult} and \eqref{al-HKR}.
Proposition \ref{prop:al-al} implies that, for every point $(n,k)$ in the 
the set $\A_m \cup \B_m$, the expression 
$$
[\al, \al] (\bsi \st^{\mo}_{n,k}) 
$$ 
may only involve  $\al(\bsi \st^{\mc}_2) = \G_{\ed}$, 
$\al(\bsi \st^{\mo}_2)= \G_{\ww}$ and the values 
$$
\al (\bsi \st^{\mo}_{n_1,k_1}) 
$$
for $(n_1, k_1) \in \A_m$\,. 

This observation allows us to formulate the following definition:
%
%
\begin{defi}
\label{dfn:m-approx}
An \emph{$m$-th approximation} to an SFQ is a degree $1$ element 
$$
\al \in \Conv(\bsi \oc, \KGra)
$$
which satisfies \eqref{al-Schouten}, \eqref{al-mult}, \eqref{al-HKR}, 
\begin{equation}
\label{al-outside}
\al(\bsi \st^{\mo}_{n,k}) = 0, \qquad \forall~~ (n,k) \notin \A_m\,,
\end{equation}
\begin{equation}
\label{MC-al-approx}
[\al, \al] (\bsi \st^{\mo}_{n,k}) = 0, \qquad \forall ~~ (n,k) \in \A_m \cup \B_m\,,
\end{equation}
and Property \ref{P:Pi-arity-2}.
\end{defi}
\begin{remark}
\label{rem:approx}
Due to \eqref{al-Schouten}, \eqref{al-mult}, and \eqref{al-mult-higher}
the equations
\begin{equation}
\label{MC-t-mc}
[\al, \al] (\bsi \st^{\mc}_{n}) = 0 \qquad \textrm{and} \qquad  [\al, \al] (\bsi \st^{\mo}_{k})  = 0
\end{equation}
hold for all $n \ge 2$ and $k \ge 2$, respectively. 
\end{remark}
\begin{remark}
\label{rem:practical}
Note that the values
$$
\al (\bsi \st^{\mo}_{n,k}) 
$$
for $(n,k) \notin \A_m$ do not show up in equation \eqref{MC-al-approx}. 
So we might as well think that the values $\al (\bsi \st^{\mo}_{n,k})$ are 
specified only for $(n,k) \in \A_m$. Let us also observe that, since 
$$
\al (\bsi \st^{\mo}_{1,k}) = \frac{1}{k!} \G^{\br}_k\,,
$$
the set 
$$
\A_m \cap \{(n,k) \in \bbZ_{\ge 2} \times \bbZ_{\ge 0}\}
$$
has finitely many points, and
$\KGra(n,k)^{\mo}$ has finite dimensional graded pieces
for every pair $(n,k)$, to introduce an $m$-th approximation $\al$ to an SFQ, 
we need to specify only finitely many coefficients in finitely 
many linear combinations 
$$
\al (\bsi \st^{\mo}_{n,k}), \qquad \textrm{for}~~~ (n,k) \in \A_m \cap \{(n,k) \in \bbZ_{\ge 2} \times \bbZ_{\ge 0}\}. 
$$
Similarly, statement 7 from Proposition \ref{prop:al-al} implies that 
\eqref{MC-al-approx} is equivalent to a finite number of polynomial equations 
on the above coefficients.  
\end{remark}

%
%

\subsubsection{There is the unique second approximation to an SFQ}
\label{sec:al-2}

Let us observe that $\A_2$ and $\B_2$ are  
$$
\A_2 = \{(1,k) ~|~ k \ge 0 \} \cup \{ (2,0), (2,1)\} 
\qquad  \textrm{and} \qquad
\B_2 = \{(3,0), (2,2)\}
$$
respectively. 

So a second approximation $\al$ to an SFQ is completely 
determined by the values: 
\begin{equation}
\label{al-2-values}
\al (\bsi \st^{\mo}_{2,0}) \in \KGra(2,0)^{\mo}\,, \qquad \al (\bsi \st^{\mo}_{2,1}) \in \KGra(2,1)^{\mo}\,.  
\end{equation}

It turns out that
\begin{prop}
\label{prop:al-2}
By setting 
\begin{equation}
\label{al-2-dfn}
\al (\bsi \st^{\mo}_{2,0}) = 0 \qquad \textrm{and} \qquad   \al (\bsi \st^{\mo}_{2,1}) = 0,
\end{equation}
we get a second approximation $\al$ to an SFQ. Moreover, if $\al$ is a 
second approximation to an SFQ then $\al$ satisfies \eqref{al-2-dfn}.
In other words, the second approximation at an SFQ is unique.  
\end{prop} 
\begin{proof}
Let $\beta^K$ be Kontsevich's SFQ from Section \ref{sec:Konts}.  
The first part of this proposition is settled by computing 
the weights of the graphs in $\beta^K (\bsi \st^{\mo}_{2,0})$ and 
$\beta^K  (\bsi \st^{\mo}_{2,1})$. This is done in Appendix \ref{app:Konts-ptys}.

According to Claim \ref{cl:Konts-2}, we have
$$
\beta^K (\bsi \st^{\mo}_{2,0})  = 0, \qquad \beta^K  (\bsi \st^{\mo}_{2,1}) = 0.
$$

Thus equations \eqref{al-2-dfn} indeed define a second approximation. 

To prove the uniqueness, we observe that 
$\al (\bsi \st^{\mo}_{2,0})$ is a degree $-2$ vector in
$$
\big( \KGra(2,0)^{\mo} \big)^{S_2}\,.
$$

It is not hard to see that the degree $-2$ component of $\KGra(2,0)^{\mo}$ is 
spanned by the graph 
$$
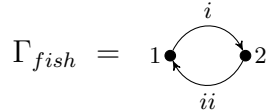
\begin{tikzpicture}[scale=0.5, >=stealth']
\tikzstyle{w}=[circle, draw, minimum size=4, inner sep=1]
\tikzstyle{b}=[circle, draw, fill, minimum size=4, inner sep=1]
\draw (-2.5,0) node[anchor=center] {{$\G_{fish}  ~ = ~ $}};
\node [b] (b1) at (0,0) {};
\draw (-0.4,0) node[anchor=center] {{\small $1$}};
\node [b] (b2) at (2,0) {};
\draw (2.4,0) node[anchor=center] {{\small $2$}};
\draw [->] (b1) ..controls (0.5,1) and (1.5,1) .. (b2);
\draw [->] (b2) ..controls (1.5,-1) and (0.5,-1) .. (b1);
\draw (1,1.2) node[anchor=center] {{\small $i$}};
\draw (1,-1.2) node[anchor=center] {{\small $ii$}};
\end{tikzpicture}
$$
and 
$$
(1,2)\, \G_{fish} = -  \G_{fish}\,.
$$
Thus 
\begin{equation}
\label{al-mo-2-0}
\al (\bsi \st^{\mo}_{2,0}) = 0 
\end{equation}
due to the fact that $\al$ is compatible with the action of $S_2$.  
  
Since multiple edges and loops are not allowed and $\al (\bsi \st^{\mo}_{2,1})$ is 
$S_2 \times \{ \id \}$-invariant, we have 
$$
\al(\bsi \st^{\mo}_{2,1})  = \la (\G_{tp} + (1,2) \G_{tp}) + \ka (\G_{\n} + (1,2) \G_{\n}),
$$
where $\G_{tp}$ (resp. $\G_{\n}$) is the first (resp. the third) graph\footnote{We 
set $1_{\mc} < 2_{\mc} < 1_{\mo}$ and use the corresponding lexicographic 
order on the set of edges.} in figure \ref{fig:mo-2-1} and $\la, \ka \in \bbQ$. 
 
Unfolding the left hand side of the identity 
$$
\al \bullet \al (\bsi \st^{\mo}_{3,0}) = 0 
$$ 
and using \eqref{cD-st-mo-nk},  \eqref{al-Schouten} and \eqref{al-mo-2-0}, we get 
\begin{equation}
\label{eq-la-kappa}
\la\, \sum_{\si \in S_3} \si (\G_{tp, 3}) ~ + ~ \ka\, \sum_{\si \in S_3} \si (\G_{\n, 3}) = 0,  
\end{equation}
where $\G_{tp, 3}$ (resp. $\G_{\n, 3}$) is the left (resp. right) graph 
in figure \ref{fig:tp-3-n-3}. 
\begin{figure}[htp] 
\centering 
\begin{minipage}[t]{0.45\linewidth}
\centering 
\begin{tikzpicture}[scale=0.5, >=stealth']
\tikzstyle{w}=[circle, draw, minimum size=4, inner sep=1]
\tikzstyle{b}=[circle, draw, fill, minimum size=4, inner sep=1]
\node [b] (b1) at (0,0) {};
\draw (-0.55,-0.15) node[anchor=center] {{\small $1$}};
\node [b] (b2) at (0,2) {};
\draw (0,2.5) node[anchor=center] {{\small $2$}};
\node [b] (b3) at (0,-1.5) {};
\draw (0,-2) node[anchor=center] {{\small $3$}};
\draw [->] (b1) ..controls (1,0.5) and (1,1.5) .. (b2);
\draw [->] (b2) ..controls (-1,1.5) and (-1, 0.5) .. (b1);
\draw [->] (b1) edge (b3);
\end{tikzpicture}
\end{minipage}
\hspace{0.1cm}
\begin{minipage}[t]{0.45\linewidth}
\centering 
\begin{tikzpicture}[scale=0.5, >=stealth']
\tikzstyle{w}=[circle, draw, minimum size=4, inner sep=1]
\tikzstyle{b}=[circle, draw, fill, minimum size=4, inner sep=1]
\node [b] (b1) at (0,0) {};
\draw (0,0.5) node[anchor=center] {{\small $1$}};
\node [b] (b2) at (2,0) {};
\draw (2,0.5) node[anchor=center] {{\small $2$}};
\node [b] (b3) at (1,-2.5) {};
\draw (1,-3) node[anchor=center] {{\small $3$}};
\draw [->] (b1) edge (b2);
\draw [->] (b1) edge (b3) (b2) edge (b3);
\end{tikzpicture}
\end{minipage}
\caption{The graphs $\G_{tp, 3}$ and $\G_{\n, 3}$} \label{fig:tp-3-n-3}
\end{figure}

Since the underlying directed graphs $\G_{tp, 3}$ and $\G_{\n, 3}$ are not isomorphic 
and both $\G_{tp, 3}$ and $\G_{\n, 3}$ have only the trivial automorphism, equation 
\eqref{eq-la-kappa} can hold only if
$$
\la = \ka  = 0. 
$$

Thus $\al (\bsi \st^{\mo}_{2,1})$ is also zero and the proposition is proved.  
\end{proof}

\subsection{The recursive construction of an SFQ} 
\label{sec:SFQ-constr} 

The main result of this paper is the following theorem which guarantees 
that a (stable) formality quasi-isomorphism for Hochschild cochains over $\bbQ$ 
can be constructed recursively:
\begin{thm}
\label{thm:main}
Let $m$ be an integer $\ge 2$ and $\al$ be an $m$-th approximation to an SFQ.
Then there exists an $(m+1)$-th approximation $\ti{\al}$ to an SFQ such that 
\begin{equation}
\label{ti-al-al}
\ti{\al} (\bsi \st^{\mo}_{n,k}) = \al (\bsi \st^{\mo}_{n,k}), \qquad \forall~~ (n,k) \in \A_m\,.
\end{equation}
Moreover, $\ti{\al}$ can be constructed by solving a finite dimensional linear system. 
In addition, a second approximation to an SFQ exists and it is unique. 
\end{thm}

The proof of this theorem is based on the following technical statement: 
\begin{prop}
\label{prop:main}
Let $m$ be an integer $\ge 2$ and $\al$ be an $m$-th approximation 
to an SFQ. Then there exists a MC element 
\begin{equation}
\label{beta-needed}
\beta \in \Conv(\bsi \oc, \KGra \otimes_{\bbQ} \bbR)
\end{equation}
satisfying \eqref{al-Schouten}, \eqref{al-mult}, \eqref{al-HKR}, 
\begin{equation}
\label{Pi-beta-2-k}
\Pi \big( \beta(\bsi \st^{\mo}_{2, k})  \big) = 0, \qquad \forall ~~ k \ge 1\,,
\end{equation}
and 
\begin{equation}
\label{beta-alpha}
\beta(\bsi \st^{\mo}_{n,k}) =  \al(\bsi \st^{\mo}_{n,k}), \qquad \forall~~ (n,k) \in \A_m\,.
\end{equation}
\end{prop}
 
We will prove this proposition in Section \ref{sec:proof-main} and now we will deduce 
Theorem \ref{thm:main} from this proposition. 

~\\

\begin{proof}[~of Theorem \ref{thm:main}] The last statement of Theorem \ref{thm:main}
follows from Proposition \ref{prop:al-2}. Thus we only need to prove the 
existence of an $(m+1)$-th approximation $\ti{\al}$ to an SFQ for 
which \eqref{ti-al-al} holds. 

For this purpose, we set $\C_m : = \A_{m+1} \setminus \A_m$, i.e.
$$
\C_m : = \{ (m+1,0), (m+1,1), (m,2), (m,3), (m-1,4), (m-1,5), \dots 
$$
\begin{equation}
\label{C-m}
 \dots ,  (2, 2m-2), (2, 2m-1) \}.
\end{equation}

Our goal is to prove that there exists a degree $1$ element $\ti{\al}$ 
in \eqref{Conv-oc-o-KGra} whose values may be different from those of $\al$
only at 
$$
\bsi \st^{\mo}_{n,k} \qquad \textrm{for}~~~ (n,k) \in \C_m\,, 
$$
and such that 
$$
\Pi \big( \ti{\al}(\bsi \st^{\mo}_{2, 2m-2})  \big) = \Pi \big( \ti{\al}(\bsi \st^{\mo}_{2, 2m-1})  \big) =  0
$$
and\footnote{It is statement 7 in Proposition \ref{prop:al-al} which implies that 
\eqref{ti-al-MC} holds automatically for $n=1$.}
\begin{equation}
\label{ti-al-MC}
[\ti{\al}, \ti{\al}](\bsi \st^{\mo}_{n,k})  = 0  \qquad  \forall ~~ (n,k) \in \C_m \cup \B_{m+1}\,.
\end{equation}

So we set 
$$
\ti{\al} (\bsi \st^{\mo}_{n,k}) : = \al (\bsi \st^{\mo}_{n,k}), \qquad \forall~~ (n,k) \in \A_m
$$
and try to find values for 
\begin{equation}
\label{unknown}
\ti{\al} (\bsi \st^{\mo}_{n,k}), \qquad \textrm{ for } ~~ (n,k) \in \C_m
\end{equation}
such that equations \eqref{ti-al-MC} hold.

Due to the statements about expressions \eqref{al-al-m-k}, \eqref{al-al-2-k}, \eqref{al-al-m-1} and 
\eqref{al-al-m-0} in Proposition \ref{prop:al-al}, the unknown vectors \eqref{unknown}
enter equations \eqref{ti-al-MC} linearly.

Since 
$$
\ti{\al} (\bsi \st^{\mo}_{n,k}) = \al (\bsi \st^{\mo}_{n,k}) =  \beta (\bsi \st^{\mo}_{n,k}) 
$$
for all $(n,k) \in \A_m$, equations \eqref{ti-al-MC} has a solution \eqref{unknown} 
(possibly over $\bbR$). Thus, since \eqref{ti-al-MC} form an inhomogeneous linear system with the 
right hand defined over $\bbQ$, a solution for this system over $\bbQ$ does exist. 

Thus $\ti{\al}$ is a desired $(m+1)$-th approximation to an SFQ
and Theorem \ref{thm:main} is proved.
\end{proof}

Theorem \ref{thm:main} allows us to produce the following algorithm. 
We start with the second approximation $\al_2$ to an SFQ presented in Section 
\ref{sec:al-2}. Due to Theorem \ref{thm:main}, there exists a third approximation 
$\al_3$ to an SFQ which extends $\al_2$ in the sense that 
$$
\al_3 (\bsi \st^{\mo}_{n,k}) = \al_2 (\bsi \st^{\mo}_{n,k}) \qquad \forall~~ (n,k) \in \A_2.
$$ 
Applying Theorem \ref{thm:main} again, we conclude that there exists 
a fourth approximation $\al_4$ which extends $\al_3$. 
At the $m$-th step of this algorithm, we take the $m$-th approximation 
to an SFQ from the previous step and produce an $(m+1)$-th approximation. 
Each step of this algorithm consists of solving a finite dimensional linear system.  

Since the union of the nested sets $\A_2 \subset \A_3 \subset \A_4 \subset \dots$ coincides with 
the whole region $\bbZ_{\ge 1} \times \bbZ_{\ge 0}$, we conclude that 
\begin{cor}
\label{cor:main}
Using the above algorithm, one can produce a stable formality 
quasi-iso\-mor\-phism for Hochschild cochains (over $\bbQ$) recursively.  
Moreover, for every $m$-th approximation $\al_m$ to an SFQ, there exists 
a MC element (corresponding at an SFQ)
$$
\al \in \Conv(\bsi \oc, \KGra)
$$
which ``extends'' $\al_m$ in the sense that 
$$
\al(\bsi \st^{\mo}_{n,k}) = \al_m(\bsi \st^{\mo}_{n,k}) \quad \forall ~ (n,k) \in \A_m\,.
$$
\end{cor}
\qed

%
%
\begin{remark}
\label{rem:globe}
Let us observe that \cite[Theorem 4.1]{exhausting} implies that any 
SFQ $\al$ (defined over $\bbQ$) is homotopy equivalent to an SFQ $\al^{\textrm{glob}}$ which 
\begin{itemize}

\item is defined over $\bbQ$, and 

\item can be used to construct a zigzag of $L_{\infty}$ quasi-isomorphisms 
connecting the sheaf of polyvector fields to the sheaf of polydifferential operators 
on an arbitrary smooth variety $X$ defined over any extension of $\bbQ$. 

\end{itemize}
This zigzag of $L_{\infty}$ quasi-isomorphisms is 
constructed using the machinery of formal geometry 
\cite{CEFT}, \cite[Sections 2,3]{Chern}, \cite{VdB-glob}, \cite{Ye}.
\end{remark}

\section{Star products modulo $(\ve^m)$ } 
\label{sec:star}

Although an $m$-th approximation $\al$ to an SFQ does not give us 
an $\infty$ morphism from the $\La\Lie$ algebra $\PV$ of polyvector 
fields to the dg $\La\Lie$ algebra $\Cbu(A)$ of Hochschild cochains
for the polynomial algebra $A$, we can still use $\al$ to construct 
approximations to star products.  

Let $\bbK$ be any field extension of $\bbQ$ and $A$ be the 
polynomial algebra 
$$
A : = \bbK[x^1, x^2, \dots, x^d]. 
$$
Let $\PV$ be the $\La\Lie$-algebra of polyvector fields on 
the corresponding affine space, i.e. 
$$
\PV  =  \bbK[x^1, x^2, \dots, x^d, \te_1, \te_2 \dots, \te_d]
$$
as the graded commutative algebra, where $\te_1, \te_2, \dots, \te_d$
are degree $1$ variables\footnote{$\PV$ should not be confused with the 
polynomial algebra in $2 d$ variables. Since $\te$'s have degree $1$, we have 
$\te_i \te_j = - \te_j \te_i$. In particular, $\te_i^2 = 0$.}.
Finally, we denote by $\ve$ a formal deformation parameter. 
 
Since the pair $(\PV, A)$ is an algebra over the operad $\KGra$, 
any $m$-th approximation of an SFQ gives us family of maps 
$$
U_{n,k} : \PV^{\otimes\, n} \otimes A^{\otimes\, k} \to A, \qquad (n,k) \in \bbZ_{\ge 1} \times \bbZ_{\ge 0}\,,
$$
\begin{equation}
\label{U-n-k}
U_{n,k}(v_1, \dots, v_n; a_1, \dots, a_k) : = 
 (-1)^{k (|v_1| + \dots + |v_n|)} 
\al(\bsi\st^{\mo}_{n,k}) \big(v_1 \otimes \dots \otimes v_n \otimes a_1 \otimes \dots \otimes a_k \big).
\end{equation}

This family assembles into a homomorphism of cocommutative coalgebras 
\begin{equation}
\label{U-from-al}
U : S( \bs^{-2}\, \PV ) \to S \big( \bs^{-2}\, \Cbu(A) \big) 
\end{equation}
such that for every $(n,k) \in \A_m \cup \B_m$, we have 
\begin{equation}
\label{U-m-relations}
\big( p \circ Q_{\Cbu(A)} \circ U  (v_1, \dots, v_n) \big) (a_1, \dots, a_k) =
\big( p \circ U  \circ Q_{\PV} (v_1, \dots, v_n) \big) (a_1, \dots, a_k),
\end{equation}
where $p$ is the projection $S\big( \bs^{-2}\, \Cbu(A) \big) \to  \bs^{-2} \Cbu(A)$, 
and $Q_{\Cbu(A)}$ (resp. $Q_{\PV}$) is the coderivation of $S\big( \bs^{-2}\, \Cbu(A) \big)$
(resp. $S( \bs^{-2}\, \PV )$) corresponding to the dg $\La\Lie$-structure on $\Cbu(A)$ 
(resp. $\PV$).  

The following theorem shows that the family of maps \eqref{U-n-k} has an interesting application. 
\begin{thm}
\label{thm:star-approx}
Let $m$ be an integer $\ge 2$,~ $\al$ be an $m$-th approximation 
to an SFQ and $\{ U_{n,k} \}_{(n,k) \in \bbZ_{\ge 1} \times \bbZ_{\ge 0}}$ be 
the above family of maps \eqref{U-n-k}. Let $\ka$ be a MC 
element\footnote{For example, $\ka = \ve \ka_1$, where $\ka_1$ is a polynomial 
Poisson structure.} of $\ve \, \PV[[\ve]]$. Then the formula 
\begin{equation}
\label{star-mod-m}
a * b : = a b + \sum_{n = 1}^{m-1} \frac{1}{n!} U_{n,2}(\underbrace{\ka, \ka, \dots, \ka}_{n~\emph{times}} ; a,b)
\end{equation}
defines an associative multiplication on $A[\ve]/(\ve^m)$. Moreover,  \eqref{star-mod-m} is 
a truncation modulo $(\ve^m)$ of an honest star product on $A[[\ve]]$. 
\end{thm}
\begin{proof}
The proof of this statement is basically an adaptation of the line of arguments in 
the proof of \cite[Proposition 2.2]{enhanced}.  

Let us recall that \eqref{star-mod-m} defines an associative multiplication on $A[\ve]/(\ve^m)$ 
if and only if the image of the element
\begin{equation}
\label{U-star-ka}
U_*(\ka) : = \sum_{n=1}^{\infty} \frac{1}{n!} p \circ U \big( (\bs^{-2}\, \ka)^n \big) \in \ve \, \Cbu(A)[[\ve]]
\end{equation}
in the quotient
\begin{equation}
\label{cL-m}
\cL_m : =  \ve \, \Cbu(A)[\ve] \,\big/\, \ve^m \, \Cbu(A)[\ve] 
\end{equation}
satisfies the MC equation.

Since $U$ is compatible with the comultiplication, we have 
\begin{equation}
\label{U-exp-ka}
U \big( \exp(\bs^{-2} \ka) -1 \big) = \exp \big( \bs^{-2} U_*(\ka) \big)  -1,
\end{equation}
where 
$$
\exp(\bs^{-2} \ka) - 1: = \sum_{n=1}^{\infty} \frac{1}{n!}  (\bs^{-2} \ka)^n
$$
is viewed as the element in the completion of $S(\bs^{-2}\, \ve \PV[[\ve]])$

Similarly, we have
\begin{equation}
\label{Q-cxp}
Q_{\PV} \big( \exp(\bs^{-2} \ka) -1 \big) = \exp(\bs^{-2} \ka)\, \bs^{-2} [\ka, \ka]_S = 0.
\end{equation}

Since \eqref{U-m-relations} holds for every $(n,k) \in \A_m \cup \B_m$ and 
terms in $ \ve^m \, \Cbu(A)[\ve] $ may be disregarded, equation \eqref{Q-cxp}
implies that 
$$
\Big( p \circ Q_{\Cbu(A)} \circ U \, \big( \exp(\bs^{-2} \ka) -1 \big) \Big) (a_1,a_2, a_3) =
$$
$$
\Big( p \circ U \circ Q_{\PV} \, \big( \exp(\bs^{-2} \ka) -1 \big) \Big)  (a_1,a_2, a_3) ~~\textrm{mod}~~ (\ve^m) ~ = ~ 0.
$$

Therefore, due to \eqref{U-exp-ka}, we have 
$$
\big( p \circ Q_{\Cbu(A)}  \big( \exp ( \bs^{-2} U_*(\ka) )  -1 \big) \big) (a_1,a_2, a_3) ~ = ~ 0 ~~\textrm{mod}~~ (\ve^m),
$$
which means that the image of $U_*(\ka)$ in the quotient \eqref{cL-m} indeed satisfies the MC equation. 

The first statement of the theorem is proved. 

The second statement follows from Corollary \ref{cor:main}. 
\end{proof}

\section{The proof of Proposition \ref{prop:main}}
\label{sec:proof-main} 

The main part of this section is devoted to the proof of the 
following statement. 
\begin{claim}
\label{cl:improve-beta}
Let $m$ and $\ti{m}$ be integers such that $m \ge \ti{m} \ge 3$, 
$\al$ be an $m$-th approximation to an SFQ and 
\begin{equation}
\label{beta-m-til}
\beta \in \Conv(\bsi \oc, \KGra \otimes_{\bbQ} \bbR)
\end{equation}
be a genuine MC element corresponding to an SFQ such that
\begin{equation}
\label{beta-al}
\beta(\bsi \st^{\mo}_{n,k}) = \al(\bsi \st^{\mo}_{n,k}) \qquad 
\forall ~~ (n,k) \in \A_{\tm - 1}. 
\end{equation}
If both $\al$ and $\beta$ satisfy Property \ref{P:Pi-arity-2}, 
then there exists a MC element (corresponding to an SFQ)
\begin{equation}
\label{beta-new}
\beta^{new} \in \Conv(\bsi \oc, \KGra \otimes_{\bbQ} \bbR)
\end{equation}
such that
\begin{equation}
\label{beta-al-better}
\beta^{new}(\bsi \st^{\mo}_{n,k}) = \al(\bsi \st^{\mo}_{n,k}) \qquad 
\forall ~~ (n,k) \in \A_{\ti{m}}. 
\end{equation}
\end{claim}
Proposition \ref{prop:main} follows easily from Proposition \ref{prop:al-2}, 
Claim \ref{cl:improve-beta} and Claim \ref{cl:Pi-arity-2} proved in Appendix \ref{app:Konts-ptys}.

\subsection{The proof of Claim \ref{cl:improve-beta}}
\label{sec:proof-main-claim}

In the course of the proof of Claim \ref{cl:improve-beta}, we will often replace a
MC element $\beta'$ in \eqref{Conv-bbR} by a one corresponding to 
a homotopy equivalent SFQ or by $\exp(\ga)\, \beta'$, where $\ga$ is 
a degree zero cocycle in $\dfGC \otimes_{\bbQ} \bbR$.
Due to Remark \ref{rem:pty-Pi-stable}, Property \ref{P:Pi-arity-2} for SFQs 
is stable both under homotopy equivalences and the action of 
the full directed graph complex. So we will tacitly assume that all 
MC elements corresponding to SFQs satisfy Property \ref{P:Pi-arity-2}.

An obvious direct computation shows that
\begin{claim}
\label{cl:MC-difference}
If $\al$ and $\beta$ are degree $1$ elements of  
$\Conv(\bsi \oc, \KGra \otimes_{\bbQ} \bbR)$ such that
$$
[\al, \al] (\bsi \st^{\mo}_{n,k}) = [\beta, \beta] (\bsi \st^{\mo}_{n,k}) = 0
$$ 
for a fixed pair $(n,k) \in \bbZ_{\ge 1} \times \bbZ_{\ge 0}$ then
\begin{equation}
\label{MC-difference}
[\al, \beta - \al ] (\bsi \st^{\mo}_{n,k}) + \frac{1}{2}[\beta - \al, \beta - \al] (\bsi \st^{\mo}_{n,k}) = 0.
\end{equation}
\qed
\end{claim}

\subsubsection{Modifying $\beta$ at $\bsi \st^{\mo}_{2,  2\ti{m} -4}$ and  $\bsi \st^{\mo}_{2,  2\ti{m} -3}$}

Due to \eqref{beta-al}, we have 
$$
[\beta - \al, \beta - \al] (\bsi \st^{\mo}_{2,  2\ti{m} -3}) =
[\beta - \al, \beta - \al] (\bsi \st^{\mo}_{2,  2\ti{m} -2}) = 0.
$$
Hence, applying \eqref{MC-difference} to the pairs $(2, 2\ti{m} -3)$ and  $(2, 2\ti{m} -2)$, 
we get 
\begin{equation}
\label{beta-al-ar-2}
[\al, \beta - \al]  (\bsi \st^{\mo}_{2,  2\ti{m} -3}) = [\al, \beta - \al]  (\bsi \st^{\mo}_{2,  2\ti{m} -2}) = 0. 
\end{equation}

Since $(\beta - \al) (\bsi \st^{\mo}_{1,k}) = 0$ for all $k$, equations in \eqref{beta-al-ar-2} are equivalent to 
\begin{equation}
\label{Hoch-ar-2}
\pa^{\Hoch}\, (\beta - \al) (\bsi \st^{\mo}_{2,  2\ti{m} -4}) = \pa^{\Hoch}\, (\beta - \al) (\bsi \st^{\mo}_{2,  2\ti{m} -3}) = 0,
\end{equation}
where $\pa^{\Hoch}$ is defined in \eqref{pa-Hoch}. 

Due to Property\footnote{Recall that Kontsevich's SFQ $\beta^K$ satisfies Property \ref{P:Pi-arity-2} due to 
Claim \ref{cl:Pi-arity-2} proved in Appendix \ref{app:Konts-ptys}.} \ref{P:Pi-arity-2} 
and \cite[Corollary A.9]{stable}, the $\pa^{\Hoch}$-cocycles
$(\beta - \al) (\bsi \st^{\mo}_{2,  2 \tm -4})$  and $(\beta - \al) (\bsi \st^{\mo}_{2,  2 \tm -3})$
are exact. In other words, there exist vectors 
\begin{equation}
\label{two-psi-s}
\psi_{2 \tm-5} \in \big( \KGra(2, 2\tm-5)^{\mo}  \otimes_{\bbQ} \bbR \big)^{S_2}\,, 
\qquad 
\psi_{2 \tm-4} \in \big( \KGra(2, 2 \tm-4)^{\mo} \otimes_{\bbQ} \bbR \big)^{S_2}
\end{equation}
of degrees $2 -2 \ti{m}$ and $1 - 2  \tm$, respectively, such that 
\begin{equation}
\label{Hoch-psi}
\pa^{\Hoch}   (\psi_{2 \tm-5}) = (\beta - \al) (\bsi \st^{\mo}_{2,  2\ti{m} -4}), 
\qquad 
\pa^{\Hoch}   (\psi_{2 \tm-4}) = (\beta - \al) (\bsi \st^{\mo}_{2,  2\ti{m} -3}). 
\end{equation}

The equations
$$
\xi_2 (\bsi \st^{\mo}_{2, 2\tm -5}) : = \psi_{2 \tm-5}\,, \qquad 
\xi'_2 (\bsi \st^{\mo}_{2, 2\tm -4}) : = \psi_{2 \tm-4}\,,  
$$
$$
\xi_2 (\bsi \st^{\mc}_{n_1}) = \xi'_2 (\bsi \st^{\mc}_{n_1})  : = 0 \quad \forall~~ n_1 \ge 2, 
$$
$$
\xi_2 (\bsi \st^{\mo}_{n, k}) : =  0 \quad \forall ~ (n,k) \neq (2, 2\tm -5), \qquad 
\xi'_2 (\bsi \st^{\mo}_{n, k}) : =  0 \quad \forall ~ (n,k) \neq (2, 2\tm -4)
$$
define degree zero vectors $\xi_2, \xi'_2 \in \Conv(\bsi \oc, \KGra \otimes_{\bbQ} \bbR)$ which satisfy 
condition \eqref{xi-cond-mc}. 

We use these vectors to produce the new MC element $\beta_2$: 
\begin{equation}
\label{beta-2}
\beta_2  : =  \exp([\xi'_2, ~]) \exp([\xi_2, ~])\, \beta.
\end{equation}
 
Due to equation \eqref{n-k1-Hoch} from Proposition \ref{prop:exp-xi-beta}
and \eqref{Hoch-psi}, we have 
\begin{equation}
\label{beta-2-good}
(\beta_2 - \al) (\bsi \st^{\mo}_{2,  2\ti{m} -4}) = (\beta_2 - \al) (\bsi \st^{\mo}_{2,  2\ti{m} -3}) = 0. 
\end{equation}
Moreover, due to the first statement of Proposition \ref{prop:exp-xi-beta}
and equations in \eqref{Hoch-psi}, 
$$
\beta_2 (\bsi \st^{\mo}_{n, k}) = \al (\bsi \st^{\mo}_{n, k}) 
$$
for all $(n,k) \in  \A_{\tm - 1}$. 

Thus 
\begin{itemize}
\item  $\beta_2$ is a MC element of $\Conv(\bsi \oc, \KGra \otimes_{\bbQ} \bbR)$ corresponding to an SFQ,

\item $\beta_2$ agrees with $\al$ at $\bsi \st^{\mo}_{n,k}$ for every  $(n,k) \in  \A_{\tm - 1}$, and

\item $\beta_2$ agrees with $\al$ at  $\bsi \st^{\mo}_{2,  2 \tm -4}$ and $\bsi \st^{\mo}_{2,  2 \tm - 3}$.  
\end{itemize}

\subsubsection{Constructing a sequence of MC elements $\{ \beta_r \}_{\, 2 \le r \le \tm-1}$ }

Let us consider the MC element $\beta_2$ as the base of our induction and assume that we 
constructed a MC element 
\begin{equation}
\label{beta-1r}
\beta_{r-1} ~\in~  \Conv(\bsi \oc, \KGra \otimes_{\bbQ} \bbR)
\end{equation} 
(corresponding to an SFQ) for $3 \le  r \le \tm  - 1$ such that 
\begin{equation}
\label{beta-1r-A1m}
\beta_{r-1} (\bsi \st^{\mo}_{n,k}) = \al (\bsi \st^{\mo}_{n,k}), \qquad \forall~~ (n,k) \in  \A_{\tm - 1}
\end{equation}
and 
\begin{equation}
\label{beta-1r-Am-almost}
\beta_{r-1} (\bsi \st^{\mo}_{n,k}) = \al (\bsi \st^{\mo}_{n,k})
\end{equation}
for $(n,k) \in  \A_{\tm}$ if $n \le r-1$. 

Our next goal is to construct a MC element $\beta_r$ in \eqref{Conv-bbR} 
corresponding to an SFQ and such that 
$$
\beta_{r} (\bsi \st^{\mo}_{n,k}) = \al (\bsi \st^{\mo}_{n,k}), \qquad \forall~~ (n,k) \in  \A_{\tm - 1}\,,
$$ 
$$
\beta_{r} (\bsi \st^{\mo}_{n,k}) = \al (\bsi \st^{\mo}_{n,k})
$$
for $(n,k) \in  \A_{\tm}$ if $n \le r$. 

Due to \eqref{beta-1r-A1m} and \eqref{beta-1r-Am-almost}, we have 
$$
[\beta_{r-1} - \al, \beta_{r-1} - \al] (\bsi \st^{\mo}_{r, 2 (\tm-r)+1})  = 0.
$$

Hence, applying \eqref{MC-difference} to the point $(r, 2 (\tm-r)+1)$, we get 
\begin{equation}
\label{r-layer}
[\al, \beta_{r-1} - \al] (\bsi \st^{\mo}_{r, 2 (\tm-r)+1}) = 0. 
\end{equation}

Since \eqref{beta-1r-Am-almost} holds for $(n,k) \in  \A_{\tm}$ if $n \le r-1$, 
equation \eqref{r-layer} is equivalent to 
$$
\pa^{\Hoch} (\beta_{r-1} - \al)  \big(\bsi \st^{\mo}_{r, 2 (\tm-r)} \big) =  0. 
$$

Therefore, due to \cite[Corollary A.9]{stable}, there exists a degree $(2- 2 \tm)$ vector 
$$
\psi_{r, 1} \in \big( \KGra(r, 2 (\tm-r)-1)^{\mo} \otimes_{\bbQ} \bbR  \big)^{S_r}
$$
such that the difference 
$$
(\beta_{r-1} - \al) (\bsi \st^{\mo}_{r,  2 (\tm-r)})- \pa^{\Hoch}   (\psi_{r,1}) 
$$
belongs to the subspace $\big( \Pi\KGra(r,  2 (\tm-r))^{\mo} \otimes_{\bbQ} \bbR \big)^{S_r}$.

It easy to see that the equations 
$$
\xi_{r,1} (\bsi \st^{\mo}_{r, 2 (\tm-r)-1}) : = \psi_{r, 1}, \quad
\xi_{r,1} (\bsi \st^{\mc}_{n_1}) : = 0 \quad \forall ~~ n_1 \ge 2,
$$
$$ 
\xi_{r,1} (\bsi \st^{\mo}_{n, k}) : = 0  \quad \forall ~~ (n,k) \neq  (r, 2 (\tm-r)-1)
$$
define a degree $0$ vector in $\Conv(\bsi \oc, \KGra \otimes_{\bbQ} \bbR)$.

Moreover, due to Proposition \ref{prop:exp-xi-beta}, the MC element 
$$
e^{[\xi_{r,1},  ~]} \beta_{r-1} 
$$

\begin{itemize}

\item agrees with $\al$ at $\bsi \st^{\mo}_{n,k}$ for every 
$(n,k) \in \A_{\tm -1}$ and for $(n,k) \in \A_{\tm}$ if $n \le r-1$,

\item the vector
$$
(e^{[\xi_{r,1}, ~]} \beta_{r-1})   \big(\bsi \st^{\mo}_{r,  2 (\tm-r)} \big)  ~ - ~ \al   (\bsi \st^{\mo}_{r,  2 (\tm-r)})
$$
belongs to the subspaces $\big( \Pi\KGra(r,  2 (\tm-r))^{\mo} \otimes_{\bbQ} \bbR \big)^{S_r}$.

\end{itemize}

Thus we may assume, without loss of generality, that the vector
$$
(\beta_{r-1} - \al)   \big(\bsi \st^{\mo}_{r,  2 (\tm-r)} \big) 
$$
belongs to the subspaces $\big( \Pi\KGra(r,  2 (\tm-r))^{\mo} \otimes_{\bbQ} \bbR \big)^{S_r}$ from the outset. 
 
Using equation \eqref{beta-1r-A1m}, the inclusion
$$
\{ (n,k) \in \bbZ_{\ge 1} \times \bbZ_{\ge 0} ~|~  2 n + k \le 2 \tm -1 \} ~\subset~ \A_{\tm -1}
$$
and the inequality $\tm \ge 3$, it is not hard to show that 
$$
[\beta_{r-1} - \al, \beta_{r-1} - \al] (\bsi \st^{\mo}_{r+1,  2 (\tm-r)-1}) = 0.
$$
 
Therefore, applying \eqref{MC-difference} to the point $(r+1,  2 (\tm-r)-1)$, we deduce that 
\begin{equation}
\label{alm-md-diff}
[\al,  \beta_{r-1} - \al] (\bsi \st^{\mo}_{r+1,  2 (\tm-r)-1}) = 0.
\end{equation}

Since $(\beta_{r-1} - \al) (\bsi \st^{\mo}_{n,k}) = 0$ for all $(n,k) \in \A_{\tm - 1}$, 
Proposition \ref{prop:al-al} implies that only the terms
$$
(\beta_{r-1} - \al) (\bsi \st^{\mo}_{r+1, 2 (\tm-r)-2}) 
\qquad \textrm{ and }  \qquad 
(\beta_{r-1} - \al) (\bsi \st^{\mo}_{r, 2 (\tm-r)})
$$
may contribute to the expression $[\al,  \beta_{r-1} - \al] (\bsi \st^{\mo}_{r+1,  2 (\tm-r)-1})$. 

More precisely, a direct computation gives us
$$
[\al,  \beta_{r-1} - \al] (\bsi \st^{\mo}_{r+1,  2 (\tm-r)-1}) =
$$
$$
\md\, (\beta_{r-1} - \al) (\bsi \st^{\mo}_{r, 2 (\tm-r)})  -  
\pa^{\Hoch}\, (\beta_{r-1} - \al) (\bsi \st^{\mo}_{r+1, 2 (\tm-r)-2}), 
$$
where $\md$ is the operator 
$$
\md : \big( \Pi \KGra(n,k)^{\mo}  \otimes_{\bbQ} \bbR \big)^{S_n} ~ \to ~  
\big( \Pi \KGra(n+1,k-1)^{\mo}  \otimes_{\bbQ} \bbR \big)^{S_{n+1}}
$$
defined by the formula\footnote{See \cite[Appendix B]{stable}.}
\begin{equation}
\label{md}
\md(\ga) =  k \, \sum_{i=1}^{n+1} \, (\tau_{n+1,i}, \id) 
\big(\, \ga  \,\c_{1, \mo}\, \G^{\br}_{0}
\, \big)\,, \qquad  \ga \in \big( \Pi\KGra(n,k)^{\mo}  \otimes_{\bbQ} \bbR \big)^{S_n}\,,
\end{equation}
$$
\tau_{n+1,i} : = (i, i+1, \dots, n, n+1).
$$

Note that, in this computation, we use the fact that 
$$
(\beta_{r-1} - \al) (\bsi \st^{\mo}_{r, 2 (\tm-r)}) \in \Pi\KGra(r,  2 (\tm-r))^{\mo} \otimes_{\bbQ} \bbR.
$$

Thus,
$$
\md\, (\beta_{r-1} - \al) (\bsi \st^{\mo}_{r, 2 (\tm-r)})  =   \pa^{\Hoch}\, (\beta_{r-1} - \al) (\bsi \st^{\mo}_{r+1, 2 (\tm-r)-2}).
$$
Hence, the second statement of \cite[Corollary A.9]{stable} implies that 
\begin{equation}
\label{md-closed-1}
\md\, (\beta_{r-1} - \al) (\bsi \st^{\mo}_{r, 2 (\tm-r)})  =   0.
\end{equation}

Since $2(\tm -r) \ge 1$, \cite[Corollary B.5]{stable} implies that there exists a vector 
$$
\vr_{1}  \in  \big( \Pi\KGra(r-1,  2 (\tm-r)+1)^{\mo} \otimes_{\bbQ} \bbR \big)^{S_{r-1}}
$$
of degree $2 - 2\tm$ such that
\begin{equation}
\label{md-vr1}
(\beta_{r-1} - \al) (\bsi \st^{\mo}_{r, 2 (\tm-r)}) = \md (\vr_{1}). 
\end{equation}

Let us now observe that the equations 
\begin{equation}
\label{ti-xi-1r-1}
\ti{\xi}_{r-1, 1} (\bsi \st^{\mo}_{r-1,  2 (\tm-r)+1}) : =  - \vr_1, \quad
\ti{\xi}_{r-1, 1} (\bsi \st^{\mc}_{n_1}) : = 0, \quad \forall ~~ n_1 \ge 2,  
\end{equation}
$$
\ti{\xi}_{r-1, 1} (\bsi \st^{\mo}_{n,k}) : = 0, \qquad \forall ~~ (n,k) \neq (r-1,  2 (\tm-r)+1) 
$$
define a degree $0$ vector in $\Conv(\bsi \oc, \KGra \otimes_{\bbQ} \bbR)$. 

Using $\ti{\xi}_{r-1, 1}$, we produce the new MC element 
\begin{equation}
\label{ti-beta-1r}
\ti{\beta}_{r-1} : =  \exp([\ti{\xi}_{r-1, 1}, ~]) \beta_{r-1}
\end{equation}
corresponding to an SFQ defined over $\bbR$. 
 
Due to Proposition \ref{prop:exp-xi-beta},  $\ti{\beta}_{r-1}$ and $\al$
still agrees at $\bsi \st^{\mo}_{n,k}$ for all $(n,k) \in \A_{\tm -1}$,
for $(n,k) \in \A_{\tm}$ if $n \le r-2$, and for $(n,k) = (r-1,  2 (\tm-r)+3)$. 

In addition, since $\vr_1$ is $\pa^{\Hoch}$-closed, equation \eqref{n-k1-Hoch}  implies that 
$$
\ti{\beta}_{r-1} (\bsi \st^{\mo}_{r-1, 2 (\tm-r)+2}) = \al (\bsi \st^{\mo}_{r-1, 2 (\tm-r)+2}).
$$ 
Thus,  $\ti{\beta}_{r-1}$ agrees with $\al$ at 
$\bsi \st^{\mo}_{n,k}$ for all $(n,k) \in \A_{\tm -1}$ and
for $(n,k) \in \A_{\tm}$ if $n \le r-1$. 

On the other hand, equations \eqref{n1-1k-broom0}, \eqref{md-vr1}
and the definition of  $\ti{\xi}_{r-1, 1}$ imply that
$$
\ti{\beta}_{r-1}(\bsi \st^{\mo}_{r, 2 (\tm-r)}) ~ = ~ \al (\bsi \st^{\mo}_{r, 2 (\tm-r)}). 
$$

Thus $\ti{\beta}_{r-1}$ agrees with $\al$ at 
$\bsi \st^{\mo}_{n,k}$ for all $(n,k) \in \A_{\tm -1}$,
for $(n,k) \in \A_{\tm}$ if $n \le r-1$, and for 
$$
(n,k) = (r, 2 (\tm-r)). 
$$

To construct the desired $\beta_r$, it remains to modify 
$\ti{\beta}_{r-1}$ so that the new MC element will also agree with $\al$
at $\bsi \st^{\mo}_{r, 2 (\tm-r)+1}$. 

Proceeding as above, we apply \eqref{MC-difference} to the point $(r, 2 (\tm-r)+2)$
and deduce that 
$$
\pa^{\Hoch} (\tbeta_{r-1} - \al)  \big(\bsi \st^{\mo}_{r, 2 (\tm-r)+1} \big) =  0. 
$$

Therefore, due to \cite[Corollary A.9]{stable}, there exists a degree $(1 - 2 \tm)$ vector 
$$
\psi_{r, 2} \in \big( \KGra(r, 2 (\tm-r))^{\mo} \otimes_{\bbQ} \bbR  \big)^{S_r}
$$
such that the difference 
$$
(\tbeta_{r-1} - \al) (\bsi \st^{\mo}_{r,  2 (\tm-r) + 1})- \pa^{\Hoch}   (\psi_{r,2}) 
$$
belongs to the subspace $\big( \Pi\KGra(r,  2 (\tm-r)+1)^{\mo} \otimes_{\bbQ} \bbR \big)^{S_r}$.

Setting
$$
\xi_{r,2} (\bsi \st^{\mo}_{r, 2 (\tm-r)}) : = \psi_{r, 2}, \quad
\xi_{r,2} (\bsi \st^{\mc}_{n_1}) : = 0 \quad \forall ~~ n_1 \ge 2,
$$
$$ 
\xi_{r,2} (\bsi \st^{\mo}_{n, k}) : = 0  \quad \forall ~~ (n,k) \neq  (r, 2 (\tm-r))
$$
we define a degree $0$ vector in $\Conv(\bsi \oc, \KGra \otimes_{\bbQ} \bbR)$.

As above, due to Proposition \ref{prop:exp-xi-beta}, the MC element 
$$
e^{[\xi_{r,2},  ~]} \tbeta_{r-1} 
$$

\begin{itemize}

\item agrees with $\al$ at $\bsi \st^{\mo}_{n,k}$ for every 
$(n,k) \in \A_{\tm -1}$, for $(n,k) \in \A_{\tm}$ if $n \le r-1$, and 
for $(n,k) = (r,  2 (\tm-r))$,

\item the vector
$$
(e^{[\xi_{r,2}, ~]} \tbeta_{r-1})   \big(\bsi \st^{\mo}_{r,  2 (\tm-r)+1} \big)  ~ - ~ \al   (\bsi \st^{\mo}_{r,  2 (\tm-r)+1})
$$
belongs to the subspaces $\big( \Pi\KGra(r,  2 (\tm-r) + 1)^{\mo} \otimes_{\bbQ} \bbR \big)^{S_r}$.

\end{itemize}

Thus we may assume, without loss of generality, that the vector
$$
(\tbeta_{r-1} - \al)   \big(\bsi \st^{\mo}_{r,  2 (\tm-r)+1} \big) 
$$
belongs to the subspaces $\big( \Pi\KGra(r,  2 (\tm-r) + 1)^{\mo} \otimes_{\bbQ} \bbR \big)^{S_r}$ from the outset. 
 
Using the equation
$$
\tbeta_{r-1} (\bsi \st^{\mo}_{n,k}) = \al (\bsi \st^{\mo}_{n,k}), \qquad \forall~~ (n,k) \in  \A_{\tm - 1}
$$
together with the inclusion
$$
\{ (n,k) \in \bbZ_{\ge 1} \times \bbZ_{\ge 0} ~|~  2 n + k \le 2 \tm -1 \} ~\subset~ \A_{\tm -1}
$$
and the inequality $\tm \ge 3$, it is not hard to show that 
$$
[\tbeta_{r-1} - \al, \tbeta_{r-1} - \al] (\bsi \st^{\mo}_{r+1,  2 (\tm-r)}) = 0.
$$
 
Therefore, applying \eqref{MC-difference} to the point $(r+1,  2 (\tm-r))$, we deduce that 
$$
[\al,  \tbeta_{r-1} - \al] (\bsi \st^{\mo}_{r+1,  2 (\tm-r)}) = 0.
$$

Then, using Proposition \ref{prop:al-al} together with the identity
$\tbeta_{r-1} (\bsi \st^{\mo}_{r, 2 (\tm-r)}) = \al (\bsi \st^{\mo}_{r, 2 (\tm-r)})$, 
we see that only the terms
$$
(\tbeta_{r-1} - \al) (\bsi \st^{\mo}_{r+1, 2 (\tm-r)-1}) 
\qquad \textrm{ and }  \qquad 
(\tbeta_{r-1} - \al) (\bsi \st^{\mo}_{r, 2 (\tm-r) + 1})
$$
may contribute to the expression $[\al,  \tbeta_{r-1} - \al] (\bsi \st^{\mo}_{r+1,  2 (\tm-r)})$. 

A direct computation gives us
$$
[\al,  \tbeta_{r-1} - \al] (\bsi \st^{\mo}_{r+1,  2 (\tm-r)}) =
$$
$$
\md\, (\tbeta_{r-1} - \al) (\bsi \st^{\mo}_{r, 2 (\tm-r)+1})  -  
\pa^{\Hoch}\, (\tbeta_{r-1} - \al) (\bsi \st^{\mo}_{r+1, 2 (\tm-r)-1}), 
$$
where $\md$ is defined in \eqref{md}.  As above, we use the fact that 
$$
(\tbeta_{r-1} - \al) (\bsi \st^{\mo}_{r, 2 (\tm-r)+1}) \in \Pi\KGra(r,  2 (\tm-r)+1)^{\mo} \otimes_{\bbQ} \bbR.
$$

Hence, as above, the second statement of \cite[Corollary A.9]{stable} implies that 
$$
\md\, (\tbeta_{r-1} - \al) (\bsi \st^{\mo}_{r, 2 (\tm-r) + 1})  =   0.
$$

Since $2(\tm -r) + 1 \ge 1$, \cite[Corollary B.5]{stable} implies that there exists a vector 
$$
\vr_{2}  \in  \big( \Pi\KGra(r-1,  2 (\tm-r)+2)^{\mo} \otimes_{\bbQ} \bbR \big)^{S_{r-1}}
$$
of degree $1 - 2m$ such that
\begin{equation}
\label{md-vr2}
(\tbeta_{r-1} - \al) (\bsi \st^{\mo}_{r, 2 (\tm-r)+1}) = \md (\vr_{2}). 
\end{equation}

Setting
\begin{equation}
\label{ti-xi-2r-1}
\ti{\xi}_{r-1, 2} (\bsi \st^{\mo}_{r-1,  2 (\tm-r)+2}) : =  - \vr_2, \quad
\ti{\xi}_{r-1, 2} (\bsi \st^{\mc}_{n_1}) : = 0, \quad \forall ~~ n_1 \ge 2,  
\end{equation}
$$
\ti{\xi}_{r-1, 2} (\bsi \st^{\mo}_{n,k}) : = 0, \qquad \forall ~~ (n,k) \neq (r-1,  2 (\tm-r)+2),
$$
we get a degree $0$ vector in $\Conv(\bsi \oc, \KGra \otimes_{\bbQ} \bbR)$. 

Using $\ti{\xi}_{r-1, 2}$, we set
\begin{equation}
\label{beta-r}
\beta_{r} : =  \exp([\ti{\xi}_{r-1, 2}, ~]) \tbeta_{r-1}
\end{equation}
and claim that this is the desired MC element. 
 
Indeed, due to Proposition \ref{prop:exp-xi-beta},  $\beta_{r}$ and $\al$
still agrees at $\bsi \st^{\mo}_{n,k}$ for all $(n,k) \in \A_{\tm -1}$,
for $(n,k) \in \A_{\tm}$ if $n \le r-2$,  for $(n,k) = (r-1,  2 (\tm-r)+2)$ 
and $(n,k) = (r,  2 (\tm-r))$. 

In addition, since $\vr_2$ is $\pa^{\Hoch}$-closed, equation \eqref{n-k1-Hoch} 
implies that 
$$
\beta_{r} (\bsi \st^{\mo}_{r-1, 2 (\tm-r)+3}) = \al (\bsi \st^{\mo}_{r-1, 2 (\tm-r)+3}).
$$ 
Therefore,  $\beta_r$ agrees with $\al$ at 
$\bsi \st^{\mo}_{n,k}$ for all $(n,k) \in \A_{\tm -1}$ and
for $(n,k) \in \A_{\tm}$ if $n \le r-1$. 

Finally, equations \eqref{n1-1k-broom0}, \eqref{md-vr2}
and the definition of  $\ti{\xi}_{r-1, 2}$ imply that
$$
\beta_r (\bsi \st^{\mo}_{r, 2 (\tm-r)+1}) ~ = ~ \al (\bsi \st^{\mo}_{r, 2 (\tm-r)+1}). 
$$

Thus $\beta_r$ is a desired MC element corresponding to an SFQ over $\bbR$ 
which agrees with $\al$ at $\bsi \st^{\mo}_{n,k}$ for all $(n,k) \in \A_{\tm -1}$,
and for $(n,k) \in \A_{\tm}$ if $n \le r$. In particular, the MC element $\beta_{\tm-1}$
agrees with $\al$ at $\bsi \st^{\mo}_{n,k}$ for all $(n,k) \in \A_{\tm}$ except 
possibly $(n,k) = (\tm,0)$ and  $(n,k) = (\tm,1)$.

\subsubsection{Getting rid of graphs with pikes in $(\tbeta- \al) \big(\bsi \st^{\mo}_{\tm, 0} \big)$} 

Let us recall that a \textit{pike} in a graph $\G \in \dgra_{n,k}$ is a black  
(i.e. $\mc$-colored) vertex of valency $1$ whose adjacent edge terminates 
at this vertex.

Setting 
$$
\tbeta  : = \beta_{\tm - 1}\,,
$$
we get a MC element in \eqref{Conv-bbR} which 
\begin{itemize}

\item corresponds to an SFQ, and 

\item agrees with $\al$ at $\bsi \st^{\mo}_{n,k}$ for 
every $(n,k) \in \A_{\tm - 1}$ and for $(n,k) \in \A_{\tm}$ if 
$n \le \tm -1$. 

\end{itemize}

In general, the linear combination
\begin{equation}
\label{de-beta-tm-0}
\tbeta(\bsi \st^{\mo}_{\tm, 0}) ~ - ~ \al (\bsi \st^{\mo}_{\tm, 0}) 
\end{equation}
may have graphs with pikes.  
So let us denote by $\de\beta^{r}_{\tm,0}$ the linear combination in 
$$
\KGra(\tm, 0)^{\mo}  \otimes_{\bbQ} \bbR
$$
which is obtained from 
the difference \eqref{de-beta-tm-0} by retaining only graphs with exactly $r$ pikes.

According to \cite[Lemma B.3]{stable}, we have 
\begin{equation}
\label{mdmdr}
\md \md^* (\de\beta^{r}_{\tm,0}) = r \de\beta^{r}_{\tm,0}\,,
\end{equation}
where the operator $\md^*$ is defined in equation (B.8) in \cite[Appendix B]{stable}. 

Thus, for the vector 
\begin{equation}
\label{xi-1tm-1}
\xi_{\tm-1,1} = - \sum_{r \ge 1} \frac{1}{r} \md^* (\de\beta^{r}_{\tm,0}) 
~ \in ~ \big( \Pi \KGra(\tm-1,1)^{\mo}  \otimes_{\bbQ} \bbR \big)^{S_{\tm-1}},
\end{equation}
the linear combination 
\begin{equation}
\label{tm-0-corrected}
\de\beta^{r}_{\tm,0} + \md (\xi_{\tm-1,1}) 
\end{equation}
does not involve graphs with pikes. 

Next, we define the degree $0$ vector 
$$
\xi \in \Conv(\bsi \oc, \KGra \otimes_{\bbQ} \bbR)
$$
by setting
\begin{equation}
\label{xi-pikes}
\xi (\bs^{-1}\,\st^{\mo}_{\tm -1,1}) : =  \xi_{\tm -1, 1}\,,
\end{equation}
$$
\xi (\bs^{-1}\,\st^{\mc}_{n_1}) : = 0, \quad  \xi (\bsi  \st^{\mo}_{k_1}) := 0, 
\quad \xi(\bsi \st^{\mo}_{n,k}) : = 0  
$$
for all $n_1, k_1 \ge 2$ and for all pairs $(n,k)$ in $\bbZ_{\ge 1} \times \bbZ_{\ge 0}$ such that 
$(n,k) \neq (\tm -1,1)$.

Then, we replace $\tbeta$ by 
\begin{equation}
\label{tbeta-pr}
\tbeta' : = \exp(\ad_{\xi}) \tbeta.
\end{equation}

Let $(n,k) \in \A_{\tm}$ with $n \le \tm - 1$. According to Proposition \ref{prop:exp-xi-beta}, 
\begin{equation}
\label{needed-4-tbeta-pr}
\tbeta'  (\bsi \st^{\mo}_{n,k}) = \al (\bsi \st^{\mo}_{n,k}) 
\end{equation}
if $(n,k) \neq (\tm-1, 2)$ and
$$
\tbeta'  (\bsi \st^{\mo}_{\tm-1 ,2}) = \al (\bsi \st^{\mo}_{\tm -1, 2}) -  \pa^{\Hoch}  \xi_{\tm -1, 1}\,. 
$$ 

On the other hand, since $\xi_{\tm -1, 1} \in \Pi \KGra(\tm-1,1)^{\mo}  \otimes_{\bbQ} \bbR$, 
$$
\pa^{\Hoch}  \xi_{\tm -1, 1} = 0 
$$
and hence \eqref{needed-4-tbeta-pr} holds for all $(n,k) \in \A_{\tm}$ with $n \le \tm -1$. 

In addition, since the linear combination \eqref{tm-0-corrected} does not 
involve graphs with pikes, equation \eqref{n1-1k-broom0} implies that 
$$
\tbeta'  (\bsi \st^{\mo}_{\tm,0}) -  \al  (\bsi \st^{\mo}_{\tm,0}) 
$$
does not involve graphs with pikes either. 

Thus we may assume, without loss of generality, that the linear combination 
\eqref{de-beta-tm-0} does not involve graphs with pikes from the outset.  
 
\subsubsection{Construction of $\beta^{new}$} 
\label{sec:constr-beta-new}

To construct the desired $\beta^{new}$ \eqref{beta-new}, we need to modify 
$\tbeta$ in such a way that the above properties of $\tbeta$ hold for 
$\beta^{new}$ and, in addition, 
$$
\beta^{new} (\bsi \st^{\mo}_{\tm, 0}) = \al (\bsi \st^{\mo}_{\tm, 0}), 
\qquad 
\beta^{new} (\bsi \st^{\mo}_{\tm, 1}) = \al (\bsi \st^{\mo}_{\tm, 1}). 
$$
 
Since $\tbeta - \al \in \cF_{\tm-1} \Conv(\bsi \oc, \KGra \otimes_{\bbQ} \bbR )$ and 
$\tm \ge 3$, we have 
\begin{equation}
\label{brack-diff-m1-0}
[\tbeta - \al, \tbeta - \al] (\bsi \st^{\mo}_{\tm+1, 0}) = 0. 
\end{equation}

Hence, using \eqref{MC-difference} and \eqref{brack-diff-m1-0}, we get 
\begin{equation}
\label{brack-al-diff}
[\al, \tbeta - \al] (\bsi \st^{\mo}_{\tm+1, 0}) = 0. 
\end{equation}

Due to the identity 
$$
\tbeta(\bsi \st^{\mo}_{n,k}) = \al (\bsi \st^{\mo}_{n,k}), \qquad \forall ~~ (n,k) \in \A_{\tm-1}\,,
$$
$$
[\al, \tbeta - \al] (\bsi \st^{\mo}_{\tm+1, 0}) ~ = ~ - \sum_{\tau \in \Sh_{2, \tm-1}} 
\tau \big(\,  (\tbeta - \al)  (\bsi \st^{\mo}_{\tm,0}) \circ_{1, \mc} \al(\st^{\mc}_2)  \,\big) 
$$
$$
+  ~  \sum_{\si \in \Sh_{1, \tm}} \si \big(\, \al(\st^{\mo}_{1,1}) \circ_{1,\mo} (\tbeta - \al)  (\bsi \st^{\mo}_{\tm,0}) \,\big)  
~ +  ~  \sum_{\si \in \Sh_{\tm, 1}} \si  \big(\,  (\tbeta - \al)  (\bsi \st^{\mo}_{\tm,1}) \circ_{1, \mo}  \al(\st^{\mo}_{1,0}) \,\big)=
$$
$$
- \sum_{\tau \in \Sh_{2, \tm-1}} 
\tau \big(\,  (\tbeta - \al)  (\bsi \st^{\mo}_{\tm,0}) \circ_{1, \mc} \G_{\ed}  \,\big) 
$$
$$
+  ~  \sum_{\si \in \Sh_{1, \tm}} \si \big(\, \G^{\br}_1 \circ_{1,\mo} (\tbeta - \al)  (\bsi \st^{\mo}_{\tm,0}) \,\big)  
~ +  ~  \sum_{\si \in \Sh_{\tm, 1}} \si  \big(\,  (\tbeta - \al)  (\bsi \st^{\mo}_{\tm,1}) \circ_{1, \mo}  \G^{\br}_0 \,\big).
$$

Thus equation \eqref{brack-al-diff} is equivalent to 
\begin{equation}
\label{eq-diff-m-0}
\sum_{\tau \in \Sh_{2, \tm-1}} 
\tau \big(\,  (\tbeta - \al)  (\bsi \st^{\mo}_{\tm,0}) \circ_{1, \mc} \G_{\ed}  \,\big) 
\end{equation}
$$
-  ~  \sum_{\si \in \Sh_{1, \tm}} \si \big(\, \G^{\br}_1 \circ_{1,\mo} (\tbeta - \al)  (\bsi \st^{\mo}_{\tm,0}) \,\big)  
~ -  ~  \sum_{\si \in \Sh_{\tm, 1}} \si  \big(\,  (\tbeta - \al)  (\bsi \st^{\mo}_{\tm,1}) \circ_{1, \mo}  \G^{\br}_0 \,\big) ~ =~  0.
$$

On other other hand, applying \eqref{MC-difference} to the point $(\tm, 2)$
and, using the above properties of $\tbeta$, we deduce that 
\begin{equation}
\label{Hoch-diff-m-2}
\pa^{\Hoch} \,  (\tbeta - \al) (\bsi \st^{\mo}_{\tm, 1}) = 0. 
\end{equation}

Hence, due to \cite[Corollary A.10]{stable}, the (only) white vertex of 
every graph in the linear combination 
$$
(\tbeta - \al) (\bsi \st^{\mo}_{\tm, 1}) 
$$ 
has valency $1$. 
Therefore, every graph in the third sum in \eqref{eq-diff-m-0} has a pike. 

Since the linear combination \eqref{de-beta-tm-0} does not 
involve graphs with pikes, 
all graphs with pikes coming from 
the first sum in \eqref{eq-diff-m-0} must cancel the third sum in \eqref{eq-diff-m-0}. 
Hence, \eqref{eq-diff-m-0} is equivalent to
\begin{equation}
\label{diff-cocycle}
[\, \G_{\ed},  (\tbeta - \al)  \big(\bsi \st^{\mo}_{\tm,0} \big)  \,] = 0, 
\end{equation}
where $(\tbeta - \al)  (\bsi \st^{\mo}_{\tm,0})$ is viewed 
as a vector in the full direct graph complex $\dfGC \otimes_{\bbQ} \bbR$ and 
$[~,~]$ is the Lie bracket on $\dfGC \otimes_{\bbQ} \bbR$ (see \cite[Section 6]{stable}). 

Since every graph in the linear combination  
$$
(\tbeta - \al)  \big( \bsi \st^{\mo}_{\tm,0} \big) 
$$
has $\tm$ vertices and $2 \tm - 2$ edges, it gives us a degree zero vector 
$$
\ga :  = (\tbeta - \al)  (\bsi \st^{\mo}_{\tm,0}) ~\in~ \cF_{\tm -1} \dfGC \otimes_{\bbQ} \bbR.
$$  
Moreover, due to \eqref{diff-cocycle}, $\ga$ is a cocycle in $\cF_{\tm -1} \dfGC  \otimes_{\bbQ} \bbR$.

Following \cite[Section 6.2]{stable}, we form the degree zero vector 
$J(\ga) \in \Conv(\bsi \oc, \KGra \otimes_{\bbQ} \bbR)$ by setting 
\begin{equation}
\label{J-ga}
J(\ga) (\bsi \st^{\mc}_n) : = \begin{cases}
 \ga  \qquad {\rm if} ~~  n  = \tm \,, \\
  0 \qquad {\rm otherwise}\,, 
\end{cases}
\end{equation}
$$
J(\ga) (\bsi \st^{\mo}_{k}) : = 0, \quad J(\ga) (\bsi \st^{\mo}_{n_1, k_1}) : = 0 
\qquad \forall ~~ k \ge 2, ~~ n_1 \ge 1, ~~ k_1 \ge 0.
$$

Let us denote by $\tbeta^{\dia}$ the new MC element 
\begin{equation}
\label{tbeta-dia}
\tbeta^{\dia} : = \exp([J(\ga),~])\, \tbeta.  
\end{equation}

Using the above defining relations of $J(\ga)$, it is easy to see that
$$
\tbeta^{\dia} (\bsi \st^{\mo}_{n,k}) = \tbeta (\bsi \st^{\mo}_{n,k}) \qquad \forall~ n \le \tm -1.
$$
Hence $\tbeta^{\dia}$ satisfies all the above properties of $\tbeta$. 

In addition,
$$
\exp([J(\ga),~])\, \tbeta (\bsi \st^{\mo}_{\tm, 0}) = \tbeta (\bsi \st^{\mo}_{\tm, 0})  + 
[J(\ga), \tbeta ] (\bsi \st^{\mo}_{\tm, 0})
$$
and
$$
[J(\ga), \tbeta ] (\bsi \st^{\mo}_{\tm, 0}) = - \G^{\br}_0 \circ_{1, \mc} \ga  = - \ga. 
$$

Thus, 
\begin{equation}
\label{tbeta-dia-good}
(\tbeta^{\dia} -\al)  \big(\bsi \st^{\mo}_{\tm, 0} \big) = 0. 
\end{equation}

In general, $(\tbeta^{\dia}  - \al) \big(\bsi \st^{\mo}_{\tm, 1} \big)$ may be non-zero. 
However, we know that 
$$
(\tbeta^{\dia} -\al)\, \big(\bsi \st^{\mo}_{\tm, 1} \big)  \in \Pi\KGra(\tm,  1)^{\mo} \otimes_{\bbQ} \bbR.
$$
In other words, the only white vertex of every graph in this linear combination is univalent.

So applying \eqref{MC-difference} to the point $(n,k) = (\tm+1, 0)$, it is easy to deduce that 
\begin{equation}
\label{md-diff-m-1}
\md \, (\tbeta^{\dia} -\al)\, \big(\bsi \st^{\mo}_{\tm, 1} \big)  ~ = ~ 0. 
\end{equation}
Therefore, due to \cite[Corollary B.5]{stable}, there exists a degree $1 - 2 \tm$ vector 
$$
\vr \in  \big(\Pi\KGra(\tm-1,  2)^{\mo} \otimes_{\bbQ} \bbR \big)^{S_{\tm-1}}
$$
such that 
\begin{equation}
\label{md-vr}
\md (\vr) = (\tbeta^{\dia} -\al)\, \big(\bsi \st^{\mo}_{\tm, 1} \big). 
\end{equation}

As above, we define a degree zero vector 
$\xi \in \Conv(\bsi \oc, \KGra \otimes_{\bbQ} \bbR)$ by setting 
$$
\xi(\bsi \st^{\mo}_{\tm-1, 2}) : = - \vr, \qquad 
\xi(\bsi \st^{\mc}_{n_1}) := 0,  \quad \forall ~~ n_1 \ge 2,
$$
$$
\xi(\bsi \st^{\mo}_{n , k}) := 0,  \qquad \forall ~~ (n,k) \neq (\tm-1, 2). 
$$

Finally, we set 
$$
\beta^{new}  : = \exp([\xi, ~])\, \tbeta^{\dia}\,.
$$

Using Proposition \ref{prop:exp-xi-beta} and the fact 
that $\vr$ is $\pa^{\Hoch}$-closed, it is easy to see that 
$\beta^{new}$ agrees with $\al$ at $\bsi \st^{\mo}_{n,k}$ for 
$(n,k) \in \A_{\tm}$ if $n \le \tm-1$ and for 
$(n,k) = (\tm, 0)$.

Furthermore, \eqref{n1-1k-broom0} and \eqref{md-vr} imply that 
$$
\beta^{new} (\bsi \st^{\mo}_{\tm, 1}) ~ = ~ \al  (\bsi \st^{\mo}_{\tm, 1}).  
$$

Thus equation \eqref{beta-al-better} holds and Claim \ref{cl:improve-beta} is proved. 

Since Claim \ref{cl:improve-beta} implies Proposition \ref{prop:main}, Theorem \ref{thm:main} also follows.

%
%

\appendix

\section{Additional properties of Kontsevich's SFQ}
\label{app:Konts-ptys}

Let us prove the following statement.
\begin{claim}
\label{cl:Konts-2}
For Kontsevich's SFQ $\beta^{K}$ from Section \ref{sec:Konts}, we have 
\begin{equation}
\label{beta-K-2}
\beta^K (\bsi \st^{\mo}_{2,0}) = 0  \qquad  \textrm{and} \qquad
\beta^K (\bsi \st^{\mo}_{2,1}) = 0. 
\end{equation}
\end{claim}
\begin{proof}
Since multiple edges and loops are not allowed, we have
$$
\beta^{K}(\bsi \st^{\mo}_{2, 0})  = W_{\G_{fish}} \G_{fish}\,,
$$
where 
$$
\begin{tikzpicture}[scale=0.5, >=stealth']
\tikzstyle{w}=[circle, draw, minimum size=4, inner sep=1]
\tikzstyle{b}=[circle, draw, fill, minimum size=4, inner sep=1]
\draw (-2.5,0) node[anchor=center] {{$\G_{fish}  ~ = ~ $}};
\node [b] (b1) at (0,0) {};
\draw (-0.4,0) node[anchor=center] {{\small $1$}};
\node [b] (b2) at (2,0) {};
\draw (2.4,0) node[anchor=center] {{\small $2$}};
\draw [->] (b1) ..controls (0.5,1) and (1.5,1) .. (b2);
\draw [->] (b2) ..controls (1.5,-1) and (0.5,-1) .. (b1);
\end{tikzpicture}
$$

According to \cite[Section 7.3.1.1]{K}, the weight\footnote{In fact, one can observe that
$\G_{fish} + (1,2) \big( \G_{fish} \big) = 0$ in $\KGra(2,0)$.} $W_{\G_{fish}} = 0$. 
So the first equation in \eqref{beta-K-2} holds. 

Again, since multiple edges and loops are not allowed, to find
$\beta^K (\bsi \st^{\mo}_{2,1})$, we need to compute the weights of the 
graphs shown in figure \ref{fig:mo-2-1}.
\begin{figure}[htp] 
\centering 
\begin{minipage}[t]{0.23\linewidth}
\centering 
\begin{tikzpicture}[scale=0.5, >=stealth']
\tikzstyle{w}=[circle, draw, minimum size=4, inner sep=1]
\tikzstyle{b}=[circle, draw, fill, minimum size=4, inner sep=1]
\node [b] (b1) at (0,0) {};
\draw (-0.55,-0.15) node[anchor=center] {{\small $1$}};
\node [b] (b2) at (0,2) {};
\draw (0,2.5) node[anchor=center] {{\small $2$}};
\node [w] (w1) at (0,-1.5) {};
\draw (0,-2) node[anchor=center] {{\small $1$}};
\draw [->] (b1) ..controls (1,0.5) and (1,1.5) .. (b2);
\draw [->] (b2) ..controls (-1,1.5) and (-1, 0.5) .. (b1);
\draw [->] (b1) edge (w1);
\end{tikzpicture}
\end{minipage}
\hspace{0.1cm}
\begin{minipage}[t]{0.23\linewidth}
\centering 
\begin{tikzpicture}[scale=0.5, >=stealth']
\tikzstyle{w}=[circle, draw, minimum size=4, inner sep=1]
\tikzstyle{b}=[circle, draw, fill, minimum size=4, inner sep=1]
\node [b] (b2) at (0,0) {};
\draw (-0.55,-0.15) node[anchor=center] {{\small $2$}};
\node [b] (b1) at (0,2) {};
\draw (0,2.5) node[anchor=center] {{\small $1$}};
\node [w] (w1) at (0,-1.5) {};
\draw (0,-2) node[anchor=center] {{\small $1$}};
\draw [->] (b2) ..controls (1,0.5) and (1,1.5) .. (b1);
\draw [->] (b1) ..controls (-1,1.5) and (-1, 0.5) .. (b2);
\draw [->] (b2) edge (w1);
\end{tikzpicture}
\end{minipage}
\hspace{0.1cm}
\begin{minipage}[t]{0.23\linewidth}
\centering 
\begin{tikzpicture}[scale=0.5, >=stealth']
\tikzstyle{w}=[circle, draw, minimum size=4, inner sep=1]
\tikzstyle{b}=[circle, draw, fill, minimum size=4, inner sep=1]
\node [b] (b1) at (0,0) {};
\draw (0,0.5) node[anchor=center] {{\small $1$}};
\node [b] (b2) at (2,0) {};
\draw (2,0.5) node[anchor=center] {{\small $2$}};
\node [w] (w1) at (1,-2.5) {};
\draw (1,-3) node[anchor=center] {{\small $1$}};
\draw [->] (b1) edge (b2);
\draw [->] (b1) edge (w1) (b2) edge (w1);
\end{tikzpicture}
\end{minipage}
\hspace{0.1cm}
\begin{minipage}[t]{0.23\linewidth}
\centering 
\begin{tikzpicture}[scale=0.5, >=stealth']
\tikzstyle{w}=[circle, draw, minimum size=4, inner sep=1]
\tikzstyle{b}=[circle, draw, fill, minimum size=4, inner sep=1]
\node [b] (b2) at (0,0) {};
\draw (0,0.5) node[anchor=center] {{\small $2$}};
\node [b] (b1) at (2,0) {};
\draw (2,0.5) node[anchor=center] {{\small $1$}};
\node [w] (w1) at (1,-2.5) {};
\draw (1,-3) node[anchor=center] {{\small $1$}};
\draw [->] (b2) edge (b1);
\draw [->] (b1) edge (w1) (b2) edge (w1);
\end{tikzpicture}
\end{minipage}
\caption{The graphs in the linear combination $\beta^K (\bsi \st^{\mo}_{2,1})$} \label{fig:mo-2-1}
\end{figure}

Since $\beta^K (\bsi \st^{\mo}_{2,1})$ is $S_2 \times \{\id\}$-invariant, it suffices to
compute the weights of the first graph and the third graph in figure \ref{fig:mo-2-1}. 

As for the first graph, we set $z_1 = \i$ and, using the argument from 
\cite[Section 7.3.1.1]{K}, we see that the corresponding weight is zero. 

Let $\G$ be the third graph in figure \ref{fig:mo-2-1} with the total order on 
the set of edges $(1_{\mc}, 2_{\mc}) < (1_{\mc}, 1_{\mo}) <  (2_{\mc}, 1_{\mo})$.  

Since every point in $C_{2,1}$ is uniquely represented by a tuple
$$
(z_1, z_2, q) \in \Conf_{2,1}, ~~~ \textrm{with}~~~ z_1 = \i, 
$$ 
the manifold $C_{2,1}$ is diffeomorphic to 
$$
\big\{ (z_2,q) ~|~ z_2 = x_2 + \i y_2  ~\in~  \bbC, ~~ q \in \bbR,  ~~ y_2 > 0, ~~z_2 \neq \i  \big\}.
$$

The formula
\begin{equation}
\label{psi-dfn}
(z_2, q) ~\mapsto~ (- \bz_2, -q)
\end{equation}
defines a diffeomorphism 
$$
\psi : C_{2,1} \to C_{2,1}
$$
and the Jacobian of this diffeomorphism with respect to 
the standard coordinates $(x_2, y_2, q)$ is $+1$. 

A direct computation shows that
$$
\psi_* \Big(\,   \bigwedge_{e \in E(\G)}  d \vf_{e}  \,\Big)  ~  = ~ - ~ \Big(\,   \bigwedge_{e \in E(\G)}  d \vf_{e}  \,\Big)
$$
and hence 
$$
\int_{C_{2,1}} \bigwedge_{e \in E(\G)}  d \vf_{e} ~ = ~ - ~  \int_{C_{2,1}} \bigwedge_{e \in E(\G)}  d \vf_{e} \,.
$$
Thus the weight of the third graph in figure \ref{fig:mo-2-1} is also zero. 

Claim \ref{cl:Konts-2} is proved. 
\end{proof}

Let us now prove that
%
%
\begin{claim}
\label{cl:Pi-arity-2}
Kontsevich's SFQ $\beta^{K}$ satisfies 
\begin{equation}
\label{Pi-beta-K-2-k}
\Pi \big( \beta^{K}(\bsi \st^{\mo}_{2, k})  \big) = 0, \qquad \forall~~ k \ge 0,
\end{equation}  
where $\Pi$ is the projection defined in \eqref{Pi}.
\end{claim}
\begin{proof}
For $k=0$ and $k = 1$, the desired statement follows readily from Claim \ref{cl:Konts-2}.

So let $\G \in \dgra_{2, k}$ with $k \ge 2$. 
If all edges terminating at white vertices of $\G$ originate from the same black vertex 
(as shown in figure \ref{fig:fish-tail}) then the argument given in \cite[Section 7.3.1.1]{K}
implies that the weight $W_{\G}$ of this graph is zero. 

\begin{figure}[htp] 
\centering 
\begin{tikzpicture}[scale=0.5, >=stealth']
\tikzstyle{w}=[circle, draw, minimum size=4, inner sep=1]
\tikzstyle{b}=[circle, draw, fill, minimum size=4, inner sep=1]
\node [b] (b1) at (0,0) {};
\draw (-0.55,-0.05) node[anchor=center] {{\small $1$}};
\node [b] (b2) at (0,2) {};
\draw (0,2.5) node[anchor=center] {{\small $2$}};
\node [w] (w1) at (-2,-3) {};
\draw (-2,-3.5) node[anchor=center] {{\small $1$}};
\node [w] (w2) at (-1,-3) {};
\draw (-1,-3.5) node[anchor=center] {{\small $2$}};
\draw (0.5,-3) node[anchor=center] {{$\dots$}};
\node [w] (wk) at (2,-3) {};
\draw (2,-3.5) node[anchor=center] {{\small $k$}};
\draw [->] (b1) ..controls (1,0.5) and (1,1.5) .. (b2);
\draw [->] (b2) ..controls (-1,1.5) and (-1, 0.5) .. (b1);
\draw [->] (b1) edge (w1) edge (w2) edge (wk);
\end{tikzpicture}
\caption{The weight of this graph is also zero} \label{fig:fish-tail}
\end{figure}
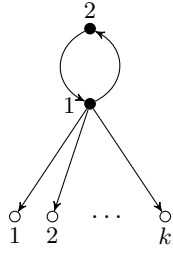

Let us now assume that each black vertex of $\G$ has valency $\ge 3$.

Recall that $C^+_{2,k}$ is the connected component of $C_{2,k}$ whose points 
are represented by tuples $(z_1, z_2, q_1, \dots, q_k)$ satisfying 
the condition
$$
q_1 < q_2 < \dots < q_k\,.
$$

We denote by  $C^{-}_{2,k}$ the connected component of $C_{2,k}$ whose points 
are represented by tuples $(z_1, z_2, q_1, \dots, q_k)$ satisfying 
the condition
$$
q_1 > q_2 > \dots > q_k\,.
$$

It is clear that the assignment 
$$
(z_1, z_2, q_1, \dots, q_k) ~\mapsto~ (- \bz_1, - \bz_2, -q_1, \dots, -q_k) 
$$
defines a diffeomorphism
\begin{equation}
\label{psi}
\psi : C^+_{2,k} \stackrel{\cong}{\longrightarrow} C^-_{2,k}
\end{equation}
whose Jacobian (with respect to the standard coordinates) is $(-1)^{k+1}$.

Let $e$ be an edge which connects the vertex corresponding to $z$ with $\Im(z) > 0$
and $q \in \bbR$. Then 
$$
d \vf_e = d \Arg(q-z) - d \Arg(q-\bz) =  d \Arg(q-z) - d \Arg(\overline{q-z}) = 2 \, d \Arg(q-z).
$$ 
For such an edge $e$, we have
$$
\psi_*  (d \vf_e ) = 2 \, d \Arg(-q + \bz) =  2 \, d \Arg(\overline{z-q}) =  -2 \, d \Arg(z-q) = - d \vf_e\,.  
$$

Moreover, since 
$$
(d \Arg(z_2 - z_1) - d \Arg(z_2- \bz_1)) \wedge (d \Arg(z_1 - z_2) - d \Arg(z_1- \bz_2)) =
$$
$$
- d \Arg(z_2 - z_1)  \wedge d \Arg(z_1- \bz_2) - d \Arg(z_2- \bz_1) \wedge d \Arg(z_1 - z_2) =
$$
$$
- 2\, d \Arg(z_2 - z_1)  \wedge d \Arg(z_1 - \bz_2) =  2\, d \Arg(z_2 - z_1)  \wedge d \Arg(z_2 - \bz_1). 
$$
and 
$$
\psi_* \big( d \Arg(z_2 - z_1)  \wedge d \Arg(z_2 - \bz_1) \big) = 
d \Arg(\bz_1 - \bz_2)  \wedge d \Arg(z_1 - \bz_2)  = 
$$
$$
(- d \Arg(z_1 - z_2)) \wedge ( - d \Arg(\bz_1 - z_2)) = d \Arg(z_2 - z_1)  \wedge d \Arg(z_2 - \bz_1),
$$
we have 
\begin{equation}
\label{psi-star-acts}
\psi_* \Big(\,   \bigwedge_{e \in E(\G)}  d \vf_{e}  \,\Big) = (-1)^{k} ~ \bigwedge_{e \in E(\G)}  d \vf_{e} \,.
\end{equation}

Combining these observations, we get\footnote{The integral over $\bC^+_{n,k}$ coincides 
with the integral over the open stratum $C^+_{n,k}$\,.} 
$$
W_{\G} = \frac{1}{(2\pi)^{k+2}} ~\int_{C^+_{2,k}}~ \bigwedge_{e \in E(\G)}  d \vf_{e} =
 \frac{(-1)^{k}}{(2\pi)^{k+2}} ~\int_{C^+_{2,k}} \psi_* \Big(\,   \bigwedge_{e \in E(\G)}  d \vf_{e}  \,\Big) =
$$
\begin{equation}
\label{W-G-psi}
- \frac{1}{(2\pi)^{k+2}} ~\int_{C^-_{2,k}}~ \bigwedge_{e \in E(\G)}  d \vf_{e}\,.
\end{equation}

On the other hand, 
\begin{equation}
\label{W-G-tilde}
\frac{1}{(2\pi)^{k+2}} ~\int_{C^-_{2,k}}~ \bigwedge_{e \in E(\G)}  d \vf_{e} ~=~
(-1)^{\frac{k(k-1)}{2}}\, W_{(\id, \si_k)\G}\,,
\end{equation}
where 
$$
\si_k = \left(
\begin{array}{ccccc}
 1 & 2  & \dots & k-1 & k   \\
 k & k-1  &  \dots & 2 & 1 
\end{array}
\right)
$$
and $(-1)^{\frac{k(k-1)}{2}}$ is precisely the sign of this permutation.

Thus \eqref{W-G-psi} and \eqref{W-G-tilde} imply that
\begin{equation}
\label{W-G-W-G-tilde}
W_{\G} = - (-1)^{\frac{k(k-1)}{2}}\, W_{(\id, \si_k)\G}\,.
\end{equation}

Let us now assume that $\G \neq (\id, \si_k)\G$ as labeled graphs. 
Then, due to \eqref{W-G-W-G-tilde}, we have 
$$
\Alt^{\mo} (W_{\G} \,\G ~ + ~  W_{(\id, \si_k)\G}  \,(\id, \si_k)\G ) =0,
$$
where $\Alt^{\mo}$ is defined in \eqref{Alt-mo}. 

Finally, if $ (\id, \si_k)\G$ coincides with\footnote{An example of a graph with 
this property is shown in figure \ref{fig:G-symm}.} 
$\G$ as the labeled graph, then, due to \eqref{W-G-W-G-tilde},
$$
\Alt^{\mo} (W_{\G} \,\G ) = 0.
$$

%
%
\begin{figure}[htp] 
\centering 
\begin{tikzpicture}[scale=0.5, >=stealth']
\tikzstyle{w}=[circle, draw, minimum size=4, inner sep=1]
\tikzstyle{b}=[circle, draw, fill, minimum size=4, inner sep=1]
\node [b] (b1) at (0,-1) {};
\draw (-0.55,-0.05) node[anchor=center] {{\small $1$}};
\node [b] (b2) at (0,2) {};
\draw (0,2.5) node[anchor=center] {{\small $2$}};
\node [w] (w1) at (-3,-3) {};
\draw (-3,-3.5) node[anchor=center] {{\small $1$}};
\node [w] (w2) at (-1,-3) {};
\draw (-1,-3.5) node[anchor=center] {{\small $2$}};
\node [w] (w3) at (1,-3) {};
\draw (1,-3.5) node[anchor=center] {{\small $3$}};
\node [w] (w4) at (3,-3) {};
\draw (3,-3.5) node[anchor=center] {{\small $4$}};
\draw [->] (b1) ..controls (0.5,0) and (0.5,1.5) .. (b2);
\draw [->] (b2) ..controls (-0.5,1.5) and (-0.5, 0) .. (b1);
\draw [->] (b1) edge (w2) edge (w3);
\draw [->] (b2) ..controls (-1.5,1.5) and (-3,-1) .. (w1);
\draw [->] (b2) ..controls (1.5,1.5) and (3,-1) .. (w4);
\end{tikzpicture}
\caption{For this graph, $(\id, \si_k)\G = \G$} \label{fig:G-symm}
\end{figure}

Claim \ref{cl:Pi-arity-2} is proved. 
\end{proof}

~\\

\noindent\textsc{Department of Mathematics,
Temple University, \\
Wachman Hall Rm. 638\\
1805 N. Broad St.,\\
Philadelphia PA, 19122 USA \\
\emph{E-mail address:} {\bf vald@temple.edu}}

\end{document}